\documentclass[12pt]{article}
\usepackage{amsfonts}
\usepackage{amssymb}

\newtheorem{theorem}{Theorem}

\begin{document}

\title{The automorphisms of class two groups of prime exponent}
\author{Michael Vaughan-Lee}
\date{January 2015}
\maketitle

\begin{abstract}
In 2012, Marcus du Sautoy and Michael Vaughan-Lee gave an example of a class
two group $G_{p}$ of prime exponent $p$ and order $p^{9}$, and they showed
that the number of descendants of $G_{p}$ of order $p^{10}$ is not a PORC
function of $p$. The fact that the number of descendants of $G_{p}$ is not
PORC is directly related to the fact that the order of the automorphism
group of $G_{p}$ is not PORC. The number of conjugacy classes of $G_{p}$ is
also not a PORC function of $p$. In this note we give a complete list of all
class two groups of prime exponent with order $p^{k}$ for $k\leq 8$. For
every group in this list we are able to show that the number of conjugacy
classes of the group is polynomial in $p$, and that the order of the
automorphism group is also polynomial in $p$. Thus, in some sense, the group 
$G_{p}$ is minimal subject to having a non-PORC number of conjugacy classes
and a non-PORC number of automorphisms.
\end{abstract}

\section{Introduction}

Graham Higman wrote two immensely important and influential papers on
enumerating $p$-groups in the late 1950s. The papers were entitled \emph{%
Enumerating }$p$\emph{-groups} I and II, and were published in the
Proceedings of the London Mathematical Society in 1960 (see \cite{higman60}
and \cite{higman60b}). In the first of these papers Higman proves that if we
let $f(p^{n})$ be the number of $p$-groups of order $p^{n}$, then $f(p^{n})$
is bounded by a polynomial in $p$. In the second of his two papers Higman
formulated his famous PORC conjecture concerning the form of the function $%
f(p^{n})$. He conjectured that for each $n$ there is an integer $N$
(depending on $n$) such that for $p$ in a fixed residue class modulo $N$ the
function $f(p^{n})$ is a polynomial in $p$. For example, {for }$p\geq 5$ the
number of groups of order $p^{6}$ is{%
\[
3p^{2}+39p+344+24\gcd (p-1,3)+11\gcd (p-1,4)+2\gcd (p-1,5). 
\]%
(See \cite{newobvl}.) So for }$p\geq 5$, $f(p^{6})$ is one of 8 polynomials
in $p$, with the choice of polynomial depending on the residue class of $p$
modulo 60. Thus $f(p^{6})$ is \textbf{P}olynomial \textbf{O}n \textbf{R}%
esidue \textbf{C}lasses. The various nineteenth century classifications of
groups of order $p^{n}$ for $n\leq 5$ show that $f(p^{n})$ is PORC for $%
n\leq 5$, and the classification of groups of order $p^{7}$ \cite{obrienvl2}
shows that $f(p^{7})$ is PORC. It is still an open question whether $%
f(p^{8}) $ is PORC, but in a recent article \cite{vlee2014} I showed that
the function enumerating the number of groups of order $p^{8}$ with exponent 
$p$ is PORC. However Marcus du Sautoy and I have found a class two group $%
G_{p}$ of order $p^{9}$ and exponent $p$ with the property that the number
of class 3 groups $H$ of order $p^{10}$ such that $H/\gamma _{3}(H)\cong
G_{p}$ is not PORC. It may still be the case that $f(p^{10})$ \emph{is}
PORC, but this example does raise a strong possibility that Higman's
conjecture fails for $n=10$. The details of this example, and a history of
the PORC conjecture can be found in \cite{dusautoyvl}. Here we give some of
the main properties of $G_{p}$.

The group $G_{p}$ is a six generator class two group of exponent $p$ ($p>2$)
with presentation%
\[
G_{p}=\left\langle 
\begin{array}{c}
x_{1},x_{2},x_{3},x_{4},x_{5},x_{6},y_{1},y_{2},y_{3}\,|\,\left[ x_{1},x_{4}%
\right] =y_{3},\,\left[ x_{1},x_{5}\right] =y_{1},\,\left[ x_{1},x_{6}\right]
=y_{2} \\ 
\left[ x_{2},x_{4}\right] =y_{1},\,\left[ x_{2},x_{5}\right] =y_{3},\,\left[
x_{3},x_{4}\right] =y_{2},\,\left[ x_{3},x_{6}\right] =y_{1}%
\end{array}%
\right\rangle 
\]%
where all other commutators are defined to be 1, and where all generators
have order $p$. The main results in \cite{dusautoyvl} are

\begin{theorem}
Let $D_{p}$ be the number of descendants of $G_{p}$ of order $p^{10}$ and
exponent $p$. Let $V_{p}$ be the number of solutions $(x,y)$ in GF$(p)$ that
satisfy $x^{4}+6x^{2}-3=0$ and $y^{2}=x^{3}-x$.

\begin{enumerate}
\item If $p=5 \bmod  12$ then $D_p=(p+1)^2/4+3$.

\item If $p=7 \bmod  12$ then $D_p=(p+1)^2/2+2$.

\item If $p=11 \bmod  12$ then $D_p=(p+1)^2/6+(p+1)/3+2$.

\item If $p=1 \bmod  12$ and $V_p=0$ then $D_p=(p+1)^2/4+3$.

\item If $p=1\bmod12$ and $V_{p}\not=0$ then $D_{p}=(p-1)^{2}/36+(p-1)/3+4$.
\end{enumerate}
\end{theorem}

\begin{theorem}
There are infinitely many primes $p=1\bmod12$ for which $V_{p}>0$. However
there is no sub-congruence of $p=1\bmod12$ for which $V_{p}>0$ for all $p$
in that sub-congruence class.
\end{theorem}

So the number of descendants of $G_{p}$ of order $p^{10}$ and exponent $p$
is not PORC, and it easily follows that the number of descendants of order $%
p^{10}$ with no restriction on the exponent is also not PORC.

We also prove in \cite{dusautoyvl} that the order of the automorphism group
of $G_{p}$ is as follows:

\begin{enumerate}
\item If $p=5\bmod12$ there are $|$GL$(2,p)|\cdot 4p^{18}$ automorphisms.

\item If $p=7\bmod12$ there are $|$GL$(2,p)|\cdot 2p^{18}$ automorphisms.

\item If $p=11\bmod12$ there are $|$GL$(2,p)|\cdot 6p^{18}$ automorphisms.

\item If $p=1\bmod12$ and $V_{p}=0$ there are $|$GL$(2,p)|\cdot 4p^{18}$
automorphisms.

\item If $p=1\bmod12$ and $V_{p}\neq 0$ there are $|$GL$(2,p)|\cdot 36p^{18}$
automorphisms.
\end{enumerate}

So the order of the automorphism group of $G_{p}$ is not PORC. To see why
the order of the automorphism group of $G_{p}$ impacts on the number of
descendants of $G_{p}$ we need to briefly recall the $p$-group generation
algorithm \cite{obrien90a}. Let $P$ be a $p$-group. The $p$-group generation
algorithm uses the lower $p$-central series, defined recursively by ${%
\mathcal{P}}_{1}(P)=P$ and ${\mathcal{P}}_{i+1}(P)=[{\mathcal{P}}_{i}(P),P]{%
\mathcal{P}}_{i}(P)^{p}$ for $i\geq 1$. The $p$-class of $P$ is the length
of this series. Each $p$-group $P$, apart from the elementary abelian ones,
is an \textit{immediate descendant} of the quotient $P/R$ where $R$ is the
last non-trivial term of the lower $p$-central series of $P$. Thus all the
groups with order $p^{n}$, except the elementary abelian one, are immediate
descendants of groups with order $p^{k}$ for $k<n$. All of the immediate
descendants of $P$ are quotients of a certain extension of $P$ (the $p$%
-covering group); the isomorphism problem for these descendants is
equivalent to the problem of determining orbits of certain subgroups of this
extension under an action of the automorphism group of $P$.

The group $G_{p}$ also has the property that the function enumerating the
number of its conjugacy classes is not PORC. I showed in \cite{vlee12} that
the number of conjugacy classes of $G_{p}$ is%
\[
p^{6}+p^{3}-1+(p^{3}-p^{2}-p+1)\times E, 
\]%
where $E$ is the number of points on the elliptic curve $y^{2}=x^{3}-x$ over
GF$(p)$ (including the point at infinity). I give a proof in \cite{vlee12}
of the well known fact that $E$ is not PORC.

It would be interesting to find more examples of finite $p$-groups with a
non-PORC number of conjugacy classes, or a non-PORC number of automorphisms,
and so I undertook a systematic search of all the class two groups of
exponent $p$ and order $p^{n}$ for $n\leq 8$. My search found nothing of
interest!

\begin{theorem}
For every prime $p>2$ there are $70$ class two groups of exponent $p$ and
order $p^{n}$ with $n\leq 8$. The number of conjugacy classes of each of
these groups is polynomial in $p$, and the number of automorphisms of each
of these groups is polynomial in $p$.
\end{theorem}

\section{The groups and their automorphisms}

For every prime $p>2$ there is one class two group of exponent $p$ of order $%
p^{3}$, one of order $p^{4}$, three of order $p^{5}$, seven of order $p^{6}$%
, fifteen of order $p^{7}$, and forty-three of order $p^{8}$. We give
presentations for all these groups below. For example the single class two
group of exponent $p$ and order $p^{3}$ has presentation%
\[
\langle a,b\,|\,\mathrm{class 2, exponent }p\rangle . 
\]%
(So it is the free group of rank 2 in the variety of class two groups of
exponent $p$.) The single class two group of exponent $p$ and order $p^{4}$
has presentation%
\[
\langle a,b,c\,|\,[c,a],\,[c,b],\mathrm{ class 2, exponent }p\rangle . 
\]

Note that the prime $p$ in these presentations is a parameter, so that the
two presentations define families of groups of order $p^{3}$ and $p^{4}$ ---
one group in each family for each prime value of $p>2$. However it is easier
to think of the presentations as defining a single group with the prime $p$
undetermined. The first of these groups has $p^{2}+p-1$ conjugacy classes,
and the second of the two groups has $p^{3}+p^{2}-p$ conjugacy classes.

Many of the presentations involve a second parameter $\omega $, which is
assumed throughout to be an integer which is primitive modulo $p$. For
example one of the class two groups of exponent $p$ and order $p^{6}$ has
presentation%
\[
\langle a,b,c,d\,|\,[c,b],\,[d,a],\,[d,b]=[c,a],\,[d,c]=[b,a]^{\omega
}\rangle . 
\]%
(For the remainder of this note we will omit the words \textquotedblleft
class 2, exponent $p$\textquotedblright\ from the presentations, taking them
as understood.) In this presentation we can take $\omega $ to be any integer
which is \emph{not} a quadratic residue modulo $p$. The isomorphism type of
the group does not depend on the particular choice of non-quadratic residue.
However, for consistency in all the different presentations it is convenient
to assume throughout that $\omega $ is always a primitive element modulo $p$.

Finally, two of the presentations involve a parameter $m$ which takes a
value such that $x^{3}+mx-1$ is irreducible over GF$(p)$ in one of the
presentations and takes a value such that $x^{3}-mx+1$ is irreducible in the
other presentation. The isomorphism types of the two groups do not depend on
the choice of $m$.

We give presentations for each of the 70 class two groups of exponent $p$
and order $p^{n}$ with $n\leq 8$, and for each of these groups we give the
polynomials giving the number of conjugacy classes and the order of the
automorphism group. We also give a description of the automorphism group
which gives the reader enough information to \textquotedblleft write
down\textquotedblright\ a set of generators for the automorphism group. If $%
G $ is a class two group of exponent $p$, and if $A$ is the automorphism
group of $G$ then $A$ induces a group of automorphisms on $G/G^{\prime }$.
If $G/G^{\prime }$ has rank $k$ then the full automorphism group of $%
G/G^{\prime }$ is GL$(k,p)$, and so the automorphisms induced by $A$ on $%
G/G^{\prime }$ form a subgroup, $B$ say, of GL$(k,p)$. Furthermore this
subgroup $B$ of GL$(k,p)$ completely determines $A$. To see this, suppose
that the defining generators of $G$ are $a_{1},a_{2},\ldots ,a_{k}$ and let $%
\alpha \in A$. If we pick arbitrary $g_{1},g_{2},\ldots ,g_{k}\in G^{\prime
} $ then there is an automorphism mapping $a_{i}$ to $(a_{i}\alpha )g_{i}$
for $i=1,2,\ldots ,k $. So $B$ completely determines $A$, and $|A|=|B|\cdot
|G^{\prime }|^{k}$.

We give three examples of the induced automorphism groups $B$. For our first
example we consider the group%
\[
\langle a,b,c,d\,|\,[c,a],\,[c,b],\,[d,a],\,[d,b],\,[d,c]\rangle 
\]%
of order $p^{5}$. The action of the automorphism group on $G/G^{\prime }$ is
given by matrices in GL$(4,p)$ of the form%
\[
\left( 
\begin{array}{cccc}
\ast & \ast & \ast & \ast \\ 
\ast & \ast & \ast & \ast \\ 
0 & 0 & \ast & \ast \\ 
0 & 0 & \ast & \ast%
\end{array}%
\right) . 
\]%
This is intended to show that the matrices in $B$ must have zeros in the
four positions shown, but that their other entries are arbitrary subject to
the restriction that the matrix must have non-zero determinant. If we write
the matrix in block form%
\[
\left( 
\begin{array}{cc}
X & Y \\ 
0 & Z%
\end{array}%
\right) 
\]%
where each block is a $2\times 2$ matrix then it is clear that $X$ and $Z$
must lie in GL$(2,p)$, whereas $Y$ is arbitrary. So the group of all these
matrices has order 
\[
(p^{2}-1)^{2}(p^{2}-p)^{2}p^{4}. 
\]
Perhaps the simplest way to generate the group is to take generators with
random entries in all the $\ast $ positions, checking that the determinant
is non-zero, and continue adding generators until they generate a subgroup
of GL$(4,p)$ of the right order. Alternatively, take generators $\left( 
\begin{array}{cc}
\omega & 0 \\ 
0 & 1%
\end{array}%
\right) $ and $\left( 
\begin{array}{cc}
-1 & 1 \\ 
-1 & 0%
\end{array}%
\right) $ for GL$(2,p)$ (where $\omega $ is a primitive element in GF$(p)$).
Then the matrix group $B$ is generated by the following matrices:%
\[
\left( 
\begin{array}{cccc}
\omega & 0 & 0 & 0 \\ 
0 & 1 & 0 & 0 \\ 
0 & 0 & 1 & 0 \\ 
0 & 0 & 0 & 1%
\end{array}%
\right) ,\,\left( 
\begin{array}{cccc}
-1 & 1 & 0 & 0 \\ 
-1 & 0 & 0 & 0 \\ 
0 & 0 & 1 & 0 \\ 
0 & 0 & 0 & 1%
\end{array}%
\right) ,\,\left( 
\begin{array}{cccc}
1 & 0 & 0 & 0 \\ 
0 & 1 & 0 & 0 \\ 
0 & 0 & \omega & 0 \\ 
0 & 0 & 0 & 1%
\end{array}%
\right) ,\,\left( 
\begin{array}{cccc}
1 & 0 & 0 & 0 \\ 
0 & 1 & 0 & 0 \\ 
0 & 0 & -1 & 1 \\ 
0 & 0 & -1 & 0%
\end{array}%
\right) , 
\]%
\[
\left( 
\begin{array}{cccc}
1 & 0 & 1 & 0 \\ 
0 & 1 & 0 & 0 \\ 
0 & 0 & 1 & 0 \\ 
0 & 0 & 0 & 1%
\end{array}%
\right) ,\,\left( 
\begin{array}{cccc}
1 & 0 & 0 & 1 \\ 
0 & 1 & 0 & 0 \\ 
0 & 0 & 1 & 0 \\ 
0 & 0 & 0 & 1%
\end{array}%
\right) ,\,\left( 
\begin{array}{cccc}
1 & 0 & 0 & 0 \\ 
0 & 1 & 1 & 0 \\ 
0 & 0 & 1 & 0 \\ 
0 & 0 & 0 & 1%
\end{array}%
\right) ,\,\left( 
\begin{array}{cccc}
1 & 0 & 0 & 0 \\ 
0 & 1 & 0 & 1 \\ 
0 & 0 & 1 & 0 \\ 
0 & 0 & 0 & 1%
\end{array}%
\right) . 
\]

As a second example, consider the group%
\[
\langle a,b,c,d\,|\,[c,b],\,[d,a],\,[d,b]=[c,a],\,[d,c]=[b,a]^{\omega
}\rangle 
\]%
of order $p^{6}$. The action of the automorphism group on $G/G^{\prime }$ is
given by matrices in GL$(4,p)$ of the form%
\[
\left( 
\begin{array}{cccc}
\alpha & \beta & \gamma & \delta \\ 
\varepsilon & \zeta & \eta & \theta \\ 
-\omega \theta & \omega \eta & \zeta & -\varepsilon \\ 
\omega \delta & -\omega \gamma & -\beta & \alpha%
\end{array}%
\right) , 
\]%
and%
\[
\left( 
\begin{array}{cccc}
\alpha & \beta & \gamma & \delta \\ 
\varepsilon & \zeta & \eta & \theta \\ 
\omega \theta & -\omega \eta & -\zeta & \varepsilon \\ 
-\omega \delta & \omega \gamma & \beta & -\alpha%
\end{array}%
\right) 
\]%
where $(\alpha ,\beta ,\gamma ,\delta )$ can be any 4-vector other than zero
($p^{4}-1$ possibilities), and where $(\varepsilon ,\zeta ,\eta ,\theta )$
can be any 4-vector which is \emph{not} in the linear span of $(\alpha
,\beta ,\gamma ,\delta )$ and $(\omega \delta ,-\omega \gamma ,-\beta
,\alpha )$ ($p^{4}-p^{2}$ possibilities). Again, one way of generating this
group of matrices is to throw in random non-singular matrices of the form
shown as generators until the required order $2(p^{4}-1)(p^{4}-p^{2})$ is
reached. Alternatively, take the matrix%
\[
\left( 
\begin{array}{cccc}
1 & 0 & 0 & 0 \\ 
0 & 1 & 0 & 0 \\ 
0 & 0 & -1 & 0 \\ 
0 & 0 & 0 & -1%
\end{array}%
\right) 
\]%
as the first generator, and then throw in random non-singular matrices of
the first form as generators until the required order is reached. A second
alternative is to first find generators for the subgroup of matrices of the
first kind with first row $(1,0,0,0)$. It is easy to see that this subgroup
is generated by matrices%
\[
\left( 
\begin{array}{cccc}
1 & 0 & 0 & 0 \\ 
1 & 1 & 0 & 0 \\ 
0 & 0 & 1 & -1 \\ 
0 & 0 & 0 & 1%
\end{array}%
\right) ,\,\left( 
\begin{array}{cccc}
1 & 0 & 0 & 0 \\ 
0 & 1 & 0 & 1 \\ 
-\omega & 0 & 1 & 0 \\ 
0 & 0 & 0 & 1%
\end{array}%
\right) ,\,\left( 
\begin{array}{cccc}
1 & 0 & 0 & 0 \\ 
0 & \zeta & \eta & 0 \\ 
0 & \omega \eta & \zeta & 0 \\ 
0 & 0 & 0 & 1%
\end{array}%
\right) 
\]%
with $\zeta ,\eta $ not both zero. Furthermore, the matrices $\left( 
\begin{array}{cc}
\zeta & \eta \\ 
\omega \eta & \zeta%
\end{array}%
\right) $ with $\zeta ,\eta $ not both zero form a group of order $p^{2}-1$
which is isomorphic to the multiplicative group of GF$(p^{2})$. So this
group is cyclic, and it is easy to find a single element which generates the
group. So we can find three matrices which generate the group of matrices of
the first type with first row $(1,0,0,0)$, and the full matrix group is then
generated by these three matrices together with 
\[
\left( 
\begin{array}{cccc}
1 & 0 & 0 & 0 \\ 
0 & 1 & 0 & 0 \\ 
0 & 0 & -1 & 0 \\ 
0 & 0 & 0 & -1%
\end{array}%
\right) 
\]%
and with matrices of the first type with general first row $(\alpha ,\beta
,\gamma ,\delta )$, and second row $(0,1,0,0)$ if $\alpha \neq 0$ or $\delta
\neq 0$, or second row $(1,0,0,0)$ if $\alpha =\delta =0$. Experimentally it
seems that we only need one of these matrices --- the one with first row $%
(0,1,0,0)$ and second row $(0,1,0,0)$.

As a third example we consider the group%
\[
\langle a,b,c,d\,|\,[c,b],\,[d,a],\,[d,b]=[c,a],\,[d,c]\rangle 
\]%
of order $p^{6}$. The action of the automorphism group on $G/G^{\prime }$ is
given by matrices in GL$(4,p)$ of the form%
\[
\left( 
\begin{array}{cccc}
\alpha & \beta & \ast & \ast \\ 
\gamma & \delta & \ast & \ast \\ 
0 & 0 & \lambda \delta & -\lambda \gamma \\ 
0 & 0 & -\lambda \beta & \lambda \alpha%
\end{array}%
\right) 
\]%
with $\lambda (\alpha \delta -\beta \gamma )\neq 0$. Here $\left( 
\begin{array}{cc}
\alpha & \beta \\ 
\gamma & \delta%
\end{array}%
\right) $ takes arbitrary values in GL$(2,p)$, $\lambda $ takes any non-zero
value, and the entries in the four positions denoted $\ast $ take arbitrary
values. This group of matrices has order $(p-1)(p^{2}-1)(p^{2}-p)p^{4}$, and
it easy to see that if we let $\omega $ be a primitive element in GF$(p)$
then it is generated by the matrices%
\[
\left( 
\begin{array}{cccc}
1 & 0 & 0 & 0 \\ 
0 & 1 & 0 & 0 \\ 
0 & 0 & \omega & 0 \\ 
0 & 0 & 0 & \omega%
\end{array}%
\right) ,\,\left( 
\begin{array}{cccc}
\omega & 0 & 0 & 0 \\ 
0 & 1 & 0 & 0 \\ 
0 & 0 & 1 & 0 \\ 
0 & 0 & 0 & \omega%
\end{array}%
\right) ,\,\left( 
\begin{array}{cccc}
-1 & 1 & 0 & 0 \\ 
-1 & 0 & 0 & 0 \\ 
0 & 0 & 0 & 1 \\ 
0 & 0 & -1 & -1%
\end{array}%
\right) , 
\]%
\[
\left( 
\begin{array}{cccc}
1 & 0 & 1 & 0 \\ 
0 & 1 & 0 & 0 \\ 
0 & 0 & 1 & 0 \\ 
0 & 0 & 0 & 1%
\end{array}%
\right) ,\,\left( 
\begin{array}{cccc}
1 & 0 & 0 & 1 \\ 
0 & 1 & 0 & 0 \\ 
0 & 0 & 1 & 0 \\ 
0 & 0 & 0 & 1%
\end{array}%
\right) ,\,\left( 
\begin{array}{cccc}
1 & 0 & 0 & 0 \\ 
0 & 1 & 1 & 0 \\ 
0 & 0 & 1 & 0 \\ 
0 & 0 & 0 & 1%
\end{array}%
\right) ,\,\left( 
\begin{array}{cccc}
1 & 0 & 0 & 0 \\ 
0 & 1 & 0 & 1 \\ 
0 & 0 & 1 & 0 \\ 
0 & 0 & 0 & 1%
\end{array}%
\right) . 
\]

The complete list of 70 class two groups of exponent $p$ and order $p^{n}$
with $n\leq 8$ follows below. For each group we give the number of conjugacy
classes, the order of the automorphism group, and a description of the
action of the automorphism group on $G/G^{\prime }$. The information shows
that in some sense the automorphism groups are independent of $p$, though of
course the entries in the matrices must lie in GF$(p)$. Also, to find sets
of generators for the matrix groups we need to make some choice of primitive
elements in GF$(p)$ and GF$(p^{2})$. There is one case, group 8.5.9, where
we need different generators for the matrix group when $p=3$, but in all
other cases the choice of generators is independent of $p$, except in the
sense just described above.

The proofs of the results below are all traditional \textquotedblleft hand
proofs\textquotedblright , albeit with machine assistance with linear
algebra over rational function fields of characteristic zero. But all the
results have been checked with a computer for small primes.
\newpage

\section{Order $p^{3}$}

\noindent\textbf{Group 3.2.1}

\[
\langle a,b\rangle 
\]

The number of conjugacy classes is $p^{2}+p-1$, and the automorphism group
has order $(p^{2}-1)(p^{2}-p)p^{2}$.

The action of the automorphism group on $G/G^{\prime }$ is given by GL$(2,p)$%
.

\section{Order $p^{4}$}

\noindent\textbf{Group 4.3.1}

\[
\langle a,b,c\,|\,[c,a],\,[c,b]\rangle 
\]

The number of conjugacy classes is $p^{3}+p^{2}-p$, and the automorphism
group has order $(p-1)(p^{2}-1)(p^{2}-p)p^{5}$.

The action of the automorphism group on $G/G^{\prime }$ is given by matrices
in GL$(3,p)$ of the form%
\[
\left( 
\begin{array}{ccc}
\ast & \ast & \ast \\ 
\ast & \ast & \ast \\ 
0 & 0 & \ast%
\end{array}%
\right) . 
\]

\section{Order $p^{5}$}

\setcounter{subsection}{2}

\subsection{Three generator groups}

\noindent\textbf{Group 5.3.1}

\[
\langle a,b,c\,|\,[c,b]\rangle 
\]

The number of conjugacy classes is $2p^{3}-p$, and the automorphism group
has order $(p-1)(p^{2}-1)(p^{2}-p)p^{8}$.

The action of the automorphism group on $G/G^{\prime }$ is given by matrices
in GL$(3,p)$ of the form%
\[
\left( 
\begin{array}{ccc}
\ast & \ast & \ast \\ 
0 & \ast & \ast \\ 
0 & \ast & \ast%
\end{array}%
\right) . 
\]

\subsection{Four generator groups}

\noindent\textbf{Group 5.4.1}

\[
\langle a,b,c,d\,|\,[c,a],\,[c,b],\,[d,a],\,[d,b],\,[d,c]\rangle 
\]

The number of conjugacy classes is $p^{4}+p^{3}-p^{2}$, and the automorphism
group has order $(p^{2}-1)^{2}(p^{2}-p)^{2}p^{8}$.

The action of the automorphism group on $G/G^{\prime }$ is given by matrices
in GL$(4,p)$ of the form%
\[
\left( 
\begin{array}{cccc}
\ast & \ast & \ast & \ast \\ 
\ast & \ast & \ast & \ast \\ 
0 & 0 & \ast & \ast \\ 
0 & 0 & \ast & \ast%
\end{array}%
\right) . 
\]

\bigskip \noindent\textbf{Group 5.4.2}

\[
\langle a,b\rangle \times _{\lbrack b,a]=[d,c]}\langle c,d\rangle 
\]

The number of conjugacy classes is $p^{4}+p-1$, and the automorphism group
has order $(p^{5}-p)(p^{5}-p^{4})(p^{3}-p)p^{2}$.

The centre of the group has order $p$. The image of $a$ can be anything
outside the centre ($p^{5}-p$ choices). The image of $b$ can be anything
outside the centralizer of the image of $a$ ($p^{5}-p^{4}$ choices). The
image of $c$ must centralize the images of $a$ and $b$ and lie outside the
centre ($p^{3}-p$ choices). The image of $d$ must centralize then images of $%
a$ and $b$, but not the image of $c$, and must be scaled so that the image
of $[d,c]$ equals the image of $[b,a]$ ($p^{2}$ choices).

\section{Order $p^{6}$}

\setcounter{subsection}{2}

\subsection{Three generator groups}

\noindent\textbf{Group 6.3.1}

\[
\langle a,b,c\rangle 
\]

The number of conjugacy classes is $p^{4}+p^{3}-p$, and the automorphism
group has order $(p^{3}-1)(p^{3}-p)(p^{3}-p^{2})p^{9}$.

The action of the automorphism group on $G/G^{\prime }$ is given by GL$(3,p)$%
.

\subsection{Four generator groups}

\noindent\textbf{Group 6.4.1}

\[
\langle a,b,c,d\,|\,[c,b],\,[d,a],\,[d,b],\,[d,c]\rangle 
\]

The number of conjugacy classes is $2p^{4}-p^{2}$ and the order of the
automorphism group is $(p-1)^{2}(p^{2}-1)(p^{2}-p)p^{13}.$

The action of the automorphism group on $G/G^{\prime }$ is given by matrices
in GL$(4,p)$ of the form%
\[
\left( 
\begin{array}{cccc}
\ast & \ast & \ast & \ast \\ 
0 & \ast & \ast & \ast \\ 
0 & \ast & \ast & \ast \\ 
0 & 0 & 0 & \ast%
\end{array}%
\right) . 
\]

\bigskip \noindent\textbf{Group 6.4.2}

\[
\langle a,b,c,d\,|\,[c,b],\,[d,a],\,[d,b]=[b,a],\,[d,c]\rangle 
\]

The number of conjugacy classes is $p^{4}+2p^{3}-p^{2}-2p+1$ and the order
of the automorphism group is $2(p^{2}-1)^{2}(p^{2}-p)^{2}p^{8}.$

The action of the automorphism group on $G/G^{\prime }$ is given by matrices
in GL$(4,p)$ of the form%
\[
\left( 
\begin{array}{cccc}
-\alpha & \beta & -\gamma & -\alpha +\delta \\ 
\varepsilon & 0 & \zeta & \varepsilon \\ 
0 & \eta & 0 & \theta \\ 
\alpha & 0 & \gamma & \alpha%
\end{array}%
\right) 
\]%
with $(\alpha \zeta -\gamma \varepsilon )(\beta \theta -\delta \eta )\neq 0$
and%
\[
\left( 
\begin{array}{cccc}
\alpha & -\beta & \gamma & \alpha -\delta \\ 
0 & \varepsilon & 0 & \zeta \\ 
\eta & 0 & \theta & \eta \\ 
0 & \beta & 0 & \delta%
\end{array}%
\right) 
\]%
with $(\alpha \theta -\gamma \eta )(\beta \zeta -\delta \varepsilon )\neq 0$.

\bigskip \noindent\textbf{Group 6.4.3}

\[
\langle a,b,c,d\,|\,[c,b],\,[d,a],\,[d,b]=[c,a],\,[d,c]\rangle 
\]

The number of conjugacy classes is $p^{4}+p^{3}-p$ and the order of the
automorphism group is $(p-1)(p^{2}-1)(p^{2}-p)p^{12}.$

The action of the automorphism group on $G/G^{\prime }$ is given by matrices
in GL$(4,p)$ of the form%
\[
\left( 
\begin{array}{cccc}
\alpha & \beta & \ast & \ast \\ 
\gamma & \delta & \ast & \ast \\ 
0 & 0 & \lambda \delta & -\lambda \gamma \\ 
0 & 0 & -\lambda \beta & \lambda \alpha%
\end{array}%
\right) 
\]%
with $\lambda (\alpha \delta -\beta \gamma )\neq 0$.

\bigskip \noindent\textbf{Group 6.4.4}

\[
\langle a,b,c,d\,|\,[c,b],\,[d,a],\,[d,b]=[c,a],\,[d,c]=[b,a]^{\omega
}\rangle 
\]

The number of conjugacy classes is $p^{4}+p^{2}-1$ and the order of the
automorphism group is $2(p^{4}-1)(p^{4}-p^{2})p^{8}.$

The action of the automorphism group on $G/G^{\prime }$ is given by matrices
in GL$(4,p)$ of the form%
\[
\left( 
\begin{array}{cccc}
\alpha & \beta & \gamma & \delta \\ 
\varepsilon & \zeta & \eta & \theta \\ 
\omega \theta & -\omega \eta & -\zeta & \varepsilon \\ 
-\omega \delta & \omega \gamma & \beta & -\alpha%
\end{array}%
\right) 
\]%
and 
\[
\left( 
\begin{array}{cccc}
\alpha & \beta & \gamma & \delta \\ 
\varepsilon & \zeta & \eta & \theta \\ 
-\omega \theta & \omega \eta & \zeta & -\varepsilon \\ 
\omega \delta & -\omega \gamma & -\beta & \alpha%
\end{array}%
\right) , 
\]%
where $(\alpha ,\beta ,\gamma ,\delta )$ can be any 4-vector other than zero
($p^{4}-1$ possibilities), and where $(\varepsilon ,\zeta ,\eta ,\theta )$
can be any 4-vector which is \emph{not} in the linear span of $(\alpha
,\beta ,\gamma ,\delta )$ and $(\omega \delta ,-\omega \gamma ,-\beta
,\alpha )$ ($p^{4}-p^{2}$ possibilities).

\subsection{Five generator groups}

\noindent\textbf{Group 6.5.1}

\[
\langle a,b\rangle \times \langle c\rangle \times \langle d\rangle \times
\langle e\rangle 
\]

The number of conjugacy classes is $p^{5}+p^{4}-p^{3}$ and the order of the
automorphism group is $%
(p^{2}-1)(p^{2}-p)(p^{3}-1)(p^{3}-p)(p^{3}-p^{2})p^{11}.$

The action of the automorphism group on $G/G^{\prime }$ is given by matrices
in GL$(5,p)$ of the form%
\[
\left( 
\begin{array}{ccccc}
\ast & \ast & \ast & \ast & \ast \\ 
\ast & \ast & \ast & \ast & \ast \\ 
0 & 0 & \ast & \ast & \ast \\ 
0 & 0 & \ast & \ast & \ast \\ 
0 & 0 & \ast & \ast & \ast%
\end{array}%
\right) . 
\]

\bigskip \noindent\textbf{Group 6.5.2}

\[
\langle a,b\rangle \times _{\lbrack b,a]=[d,c]}\langle c,d\rangle \times
\langle e\rangle 
\]

The number of conjugacy classes is $p^{5}+p^{2}-p$ and the order of the
automorphism group is $%
(p^{6}-p^{2})(p^{6}-p^{5})(p^{4}-p^{2})p^{3}(p^{2}-p). $

The centre of the group has order $p^{2}$. The image of $a$ can be anything
not in the centre ($p^{6}-p^{2}$ choices). The image of $b$ can be anything
which does not centralize the image of $a$ ($p^{6}-p^{5}$ choices). The
image of $c$ must centralize the images of $a$ and $b$ but lie outside the
centre ($p^{4}-p^{2}$ choices). The image of $d$ must centralize the images
of $a$ and $b$, but not the image of $c$, and must be scaled so that the
image of $[d,c]$ equals the image of $[b,a]$ ($p^{3}$ choices). The image of 
$e$ must lie in the centre, but not in the derived group ($p^{2}-p$ choices).

\section{Order $p^{7}$}

\setcounter{subsection}{3}

\subsection{Four generator groups}

\noindent\textbf{Group 7.4.1}

\[
\langle a,b,c,d\,|\,[c,b],\,[d,b],\,[d,c]\rangle 
\]

The number of conjugacy classes is $p^{5}+p^{4}-p^{2}$, and the automorphism
group has order $(p-1)(p^{3}-1)(p^{3}-p)(p^{3}-p^{2})p^{15}$.

The action of the automorphism group on $G/G^{\prime }$ is given by matrices
in GL$(4,p)$ of the form%
\[
\left( 
\begin{array}{cccc}
\ast & \ast & \ast & \ast \\ 
0 & \ast & \ast & \ast \\ 
0 & \ast & \ast & \ast \\ 
0 & \ast & \ast & \ast%
\end{array}%
\right) . 
\]

\newpage \noindent\textbf{Group 7.4.2}

\[
\langle a,b,c,d\,|\,[d,a],\,[d,b],\,[d,c]\rangle 
\]

The number of conjugacy classes is $p^{5}+p^{4}-p^{2}$, and the automorphism
group has order $(p-1)(p^{3}-1)(p^{3}-p)(p^{3}-p^{2})p^{15}$.

The action of the automorphism group on $G/G^{\prime }$ is given by matrices
in GL$(4,p)$ of the form%
\[
\left( 
\begin{array}{cccc}
\ast & \ast & \ast & \ast \\ 
\ast & \ast & \ast & \ast \\ 
\ast & \ast & \ast & \ast \\ 
0 & 0 & 0 & \ast%
\end{array}%
\right) . 
\]

\bigskip \noindent\textbf{Group 7.4.3}

\[
\langle a,b,c,d\,|\,[d,a],\,[c,b],\,[d,c]\rangle 
\]

The number of conjugacy classes is $3p^{4}-3p^{2}+p$, and the automorphism
group has order $2(p-1)^{4}p^{16}$.

The action of the automorphism group on $G/G^{\prime }$ is given by matrices
in GL$(4,p)$ of the form%
\[
\left( 
\begin{array}{cccc}
\ast & 0 & \ast & \ast \\ 
0 & \ast & \ast & \ast \\ 
0 & 0 & \ast & 0 \\ 
0 & 0 & 0 & \ast%
\end{array}%
\right) \mathrm{ and }\left( 
\begin{array}{cccc}
0 & \ast & \ast & \ast \\ 
\ast & 0 & \ast & \ast \\ 
0 & 0 & 0 & \ast \\ 
0 & 0 & \ast & 0%
\end{array}%
\right) . 
\]

\bigskip \noindent\textbf{Group 7.4.4}

\[
\langle a,b,c,d\,|\,[d,a]=[c,b],\,[d,b],\,[d,c]\rangle 
\]

The number of conjugacy classes is $2p^{4}+p^{3}-2p^{2}$, and the
automorphism group has order $(p-1)(p^{2}-1)(p^{2}-p)p^{17}$.

The action of the automorphism group on $G/G^{\prime }$ is given by matrices
in GL$(4,p)$ of the form%
\[
\left( 
\begin{array}{cccc}
\lambda & \ast & \ast & \ast \\ 
0 & \alpha & \beta & \ast \\ 
0 & \gamma & \delta & \ast \\ 
0 & 0 & 0 & \lambda ^{-1}(\alpha \delta -\beta \gamma )%
\end{array}%
\right) , 
\]%
with $\lambda (\alpha \delta -\beta \gamma )\neq 0$.

\bigskip \noindent\textbf{Group 7.4.5}

\[
\langle a,b,c,d\,|\,[c,a],\,[d,a]=[c,b],\,[d,b]\rangle 
\]

The number of conjugacy classes is $2p^{4}+p^{3}-2p^{2}$, and the
automorphism group has order $(p^{2}-1)^{2}(p^{2}-p)p^{13}$.

The action of the automorphism group on $G/G^{\prime }$ is given by matrices
in GL$(4,p)$ of the form%
\[
\left( 
\begin{array}{cccc}
\lambda \alpha & \lambda \beta & \mu \alpha & -\mu \beta \\ 
-\lambda \gamma & -\lambda \delta & -\mu \gamma & \mu \delta \\ 
\nu \alpha & \nu \beta & \xi \alpha & -\xi \beta \\ 
\nu \gamma & \nu \delta & \xi \gamma & -\xi \delta%
\end{array}%
\right) , 
\]%
with $(\alpha \delta -\beta \gamma )(\lambda \xi -\mu \nu )\neq 0$

\bigskip \noindent\textbf{Group 7.4.6}

\[
\langle a,b,c,d\,|\,[d,b]=[c,a]^{\omega },\,[d,c]=[b,a],\,[c,b]\rangle 
\]

The number of conjugacy classes is $p^{4}+2p^{3}-p^{2}-p$, and the
automorphism group has order $2(p^{2}-1)^{2}p^{16}$.

The action of the automorphism group on $G/G^{\prime }$ is given by matrices
in GL$(4,p)$ of the form%
\[
\left( 
\begin{array}{cccc}
\alpha & \ast & \ast & \beta \\ 
0 & \gamma & \omega \delta & 0 \\ 
0 & \delta & \gamma & 0 \\ 
\omega \beta & \ast & \ast & \alpha%
\end{array}%
\right) \mathrm{ and }\left( 
\begin{array}{cccc}
\alpha & \ast & \ast & \beta \\ 
0 & \gamma & \omega \delta & 0 \\ 
0 & -\delta & -\gamma & 0 \\ 
-\omega \beta & \ast & \ast & -\alpha%
\end{array}%
\right) , 
\]%
with $\alpha $ and $\beta $ not both zero, and with $\gamma $ and $\delta $
not both zero.

\subsection{Five generator groups}

\noindent\textbf{Group 7.5.1}

\[
\langle
a,b,c,d,e\,|\,[c,b],\,[d,a],\,[d,b],\,[d,c],\,[e,a],\,[e,b],\,[e,c],\,[e,d]%
\rangle 
\]

The number of conjugacy classes is $2p^{5}-p^{3}$ and the order of the
automorphism group is $(p-1)(p^{2}-1)^{2}(p^{2}-p)^{2}p^{18}.$

The action of the automorphism group on $G/G^{\prime }$ is given by matrices
in GL$(5,p)$ of the form%
\[
\left( 
\begin{array}{ccccc}
\ast & \ast & \ast & \ast & \ast \\ 
0 & \ast & \ast & \ast & \ast \\ 
0 & \ast & \ast & \ast & \ast \\ 
0 & 0 & 0 & \ast & \ast \\ 
0 & 0 & 0 & \ast & \ast%
\end{array}%
\right) . 
\]

\bigskip \noindent\textbf{Group 7.5.2}

\[
\langle
a,b,c,d,e\,|\,[c,b],\,[d,a],\,[d,b]=[b,a],\,[d,c],\,[e,a],\,[e,b],\,[e,c],%
\,[e,d]\rangle 
\]

The number of conjugacy classes is $p^{5}+2p^{4}-p^{3}-2p^{2}+p$ and the
order of the automorphism group is $2(p-1)(p^{2}-1)^{2}(p^{2}-p)^{2}p^{14}.$

The action of the automorphism group on $G/G^{\prime }$ is given by matrices
in GL$(5,p)$ of the form%
\[
\left( 
\begin{array}{ccccc}
-\alpha & \beta & -\gamma & -\alpha +\delta & \ast \\ 
\varepsilon & 0 & \zeta & \varepsilon & \ast \\ 
0 & \eta & 0 & \theta & \ast \\ 
\alpha & 0 & \gamma & \alpha & \ast \\ 
0 & 0 & 0 & 0 & \ast%
\end{array}%
\right) 
\]%
with $(\alpha \zeta -\gamma \varepsilon )(\beta \theta -\delta \eta )\neq 0$
and%
\[
\left( 
\begin{array}{ccccc}
\alpha & -\beta & \gamma & \alpha -\delta & \ast \\ 
0 & \varepsilon & 0 & \zeta & \ast \\ 
\eta & 0 & \theta & \eta & \ast \\ 
0 & \beta & 0 & \delta & \ast \\ 
0 & 0 & 0 & 0 & \ast%
\end{array}%
\right) 
\]%
with $(\alpha \theta -\gamma \eta )(\beta \zeta -\delta \varepsilon )\neq 0$.

\bigskip \noindent\textbf{Group 7.5.3}

\[
\langle
a,b,c,d,e\,|\,[c,b],\,[d,a],\,[d,b]=[c,a],\,[d,c],\,[e,a],\,[e,b],\,[e,c],%
\,[e,d]\rangle 
\]

The number of conjugacy classes is $p^{5}+p^{4}-p^{2}$ and the order of the
automorphism group is $(p-1)^{2}(p^{2}-1)(p^{2}-p)p^{18}.$

The action of the automorphism group on $G/G^{\prime }$ is given by matrices
in GL$(5,p)$ of the form%
\[
\left( 
\begin{array}{ccccc}
\alpha & \beta & \ast & \ast & \ast \\ 
\gamma & \delta & \ast & \ast & \ast \\ 
0 & 0 & \lambda \delta & -\lambda \gamma & \ast \\ 
0 & 0 & -\lambda \beta & \lambda \alpha & \ast \\ 
0 & 0 & 0 & 0 & \ast%
\end{array}%
\right) . 
\]%
with $\lambda (\alpha \delta -\beta \gamma )\neq 0$.

\bigskip \noindent\textbf{Group 7.5.4}

\[
\langle a,b,c,d,e\,|\,[c,b],\,[d,a],\,[d,b]=[c,a],\,[d,c]=[b,a]^{\omega
},\,[e,a],\,[e,b],\,[e,c],\,[e,d]\rangle 
\]

The number of conjugacy classes is $p^{5}+p^{3}-p$ and the order of the
automorphism group is $2(p-1)(p^{4}-1)(p^{4}-p^{2})p^{14}.$

The action of the automorphism group on $G/G^{\prime }$ is given by matrices
in GL$(5,p)$ of the form%
\[
\left( 
\begin{array}{ccccc}
\alpha & \beta & \gamma & \delta & \ast \\ 
\varepsilon & \zeta & \eta & \theta & \ast \\ 
\omega \theta & -\omega \eta & -\zeta & \varepsilon & \ast \\ 
-\omega \delta & \omega \gamma & \beta & -\alpha & \ast \\ 
0 & 0 & 0 & 0 & \ast%
\end{array}%
\right) 
\]%
and 
\[
\left( 
\begin{array}{ccccc}
\alpha & \beta & \gamma & \delta & \ast \\ 
\varepsilon & \zeta & \eta & \theta & \ast \\ 
-\omega \theta & \omega \eta & \zeta & -\varepsilon & \ast \\ 
\omega \delta & -\omega \gamma & -\beta & \alpha & \ast \\ 
0 & 0 & 0 & 0 & \ast%
\end{array}%
\right) , 
\]%
where $(\alpha ,\beta ,\gamma ,\delta )$ can be any 4-vector other than zero
($p^{4}-1$ possibilities), and where $(\varepsilon ,\zeta ,\eta ,\theta )$
can be any 4-vector which is \emph{not} in the linear span of $(\alpha
,\beta ,\gamma ,\delta )$ and $(\omega \delta ,-\omega \gamma ,-\beta
,\alpha )$ ($p^{4}-p^{2}$ possibilities).

\bigskip \noindent\textbf{Group 7.5.5}

\[
\langle
a,b,c,d,e\,|\,[c,b],\,[d,a],\,[d,b],\,[d,c],\,[e,a],\,[e,b],\,[e,c],%
\,[e,d]=[b,a]\rangle 
\]

The number of conjugacy classes is $p^{5}+p^{4}-p^{2}$ and the order of the
automorphism group is $(p-1)^{2}(p^{2}-1)(p^{2}-p)p^{15}$.

The action of the automorphism group on $G/G^{\prime }$ is given by matrices
in GL$(5,p)$ of the form%
\[
\left( 
\begin{array}{ccccc}
\alpha & \beta & \gamma & \delta & \varepsilon \\ 
0 & \alpha ^{-1}(\lambda \xi -\mu \nu ) & 0 & 0 & 0 \\ 
0 & \zeta & \eta & 0 & 0 \\ 
0 & \alpha ^{-1}(-\delta \mu +\varepsilon \lambda ) & 0 & \lambda & \mu \\ 
0 & \alpha ^{-1}(-\delta \xi +\varepsilon \nu ) & 0 & \nu & \xi%
\end{array}%
\right) 
\]%
with $\alpha ,\eta ,\lambda \xi -\mu \nu \neq 0$.

\bigskip \noindent\textbf{Group 7.5.6}

\[
\langle
a,b,c,d,e\,|\,[c,b],\,[d,a],\,[d,b]=[c,a],\,[d,c],\,[e,a],\,[e,b],\,[e,c],%
\,[e,d]=[b,a]\rangle 
\]

The number of conjugacy classes is $p^{5}+p^{3}-p$ and the order of the
automorphism group is $(p-1)(p^{2}-1)(p^{2}-p)p^{14}$.

The action of the automorphism group on $G/G^{\prime }$ is given by a
subgroup $H$ of GL$(5,p)$, where the matrices in $H$ have the form%
\[
\left( 
\begin{array}{ccccc}
\alpha & \ast & \ast & \beta & \ast \\ 
0 & \ast & \ast & 0 & \ast \\ 
0 & \ast & \ast & 0 & \ast \\ 
\gamma & \ast & \ast & \delta & \ast \\ 
0 & \ast & \ast & 0 & \ast%
\end{array}%
\right) . 
\]%
As we run through the elements of $H$, $\left( 
\begin{array}{cc}
\alpha & \beta \\ 
\gamma & \delta%
\end{array}%
\right) $ takes all values in GL$(2,p)$. If we take generators $\left( 
\begin{array}{cc}
\omega & 0 \\ 
0 & 1%
\end{array}%
\right) $ and $\left( 
\begin{array}{cc}
-1 & 1 \\ 
-1 & 0%
\end{array}%
\right) $ for GL$(2,p)$ then we obtain the following generating matrices for 
$H$:%
\[
\left( 
\begin{array}{ccccc}
\omega & \alpha & \beta & 0 & \gamma \\ 
0 & \omega \lambda & 0 & 0 & 0 \\ 
0 & 0 & \lambda & 0 & 0 \\ 
0 & \omega ^{-1}\gamma & -\omega ^{-1}\alpha & 1 & \delta \\ 
0 & 0 & 0 & 0 & \omega ^{2}\lambda%
\end{array}%
\right) \;(\lambda \neq 0) 
\]%
and%
\[
\left( 
\begin{array}{ccccc}
-1 & \alpha & \beta & 1 & \gamma \\ 
0 & \lambda & 0 & 0 & 2\lambda \\ 
0 & 0 & 0 & 0 & \lambda \\ 
-1 & \delta & \beta -\delta & 0 & \delta -\alpha \\ 
0 & -\lambda & \lambda & 0 & -\lambda%
\end{array}%
\right) \;(\lambda \neq 0). 
\]

\subsection{Six generator groups}

\noindent\textbf{Group 7.6.1}

\[
\langle a,b\rangle \times \langle c\rangle \times \langle d\rangle \times
\langle e\rangle \times \langle f\rangle 
\]

The number of conjugacy classes is $p^{6}+p^{5}-p^{4}$, and the automorphism
group has order $%
(p^{7}-p^{5})(p^{7}-p^{6})(p^{5}-p)(p^{5}-p^{2})(p^{5}-p^{3})(p^{5}-p^{4})$.

\bigskip \noindent\textbf{Group 7.6.2}

\[
\langle a,b\rangle \times _{\lbrack b,a]=[d,c]}\langle c,d\rangle \times
\langle e\rangle \times \langle f\rangle 
\]

The number of conjugacy classes is $p^{6}+p^{3}-p^{2}$, and the automorphism
group has order $%
(p^{7}-p^{3})(p^{7}-p^{6})(p^{5}-p^{3})p^{4}(p^{3}-p)(p^{3}-p^{2})$.

\bigskip \noindent\textbf{Group 7.6.3}

\[
\langle a,b\rangle \times _{\lbrack b,a]=[d,c]=[f,e]}\langle c,d\rangle
\times _{\lbrack b,a]=[d,c]=[f,e]}\langle e,f\rangle 
\]

The number of conjugacy classes is $p^{6}+p-1$, and the automorphism group
has order $(p^{7}-p)(p^{7}-p^{6})(p^{5}-p)p^{4}(p^{3}-p)p^{2}$.

\section{Order $p^{8}$}

\setcounter{subsection}{3}

\subsection{Four generator groups}

\noindent\textbf{Group 8.4.1}

\[
\langle a,b,c,d\,|\,[b,a],\,[c,a]\rangle 
\]

The number of conjugacy classes is $2p^{5}+p^{4}-2p^{3}$, and the
automorphism group has order $(p-1)^{2}(p^{2}-1)(p^{2}-p)p^{21}$.

The action of the automorphism group on $G/G^{\prime }$ is given by matrices
in GL$(4,p)$ of the form%
\[
\left( 
\begin{array}{cccc}
\ast & 0 & 0 & 0 \\ 
\ast & \ast & \ast & 0 \\ 
\ast & \ast & \ast & 0 \\ 
\ast & \ast & \ast & \ast%
\end{array}%
\right) . 
\]

\bigskip \noindent\textbf{Group 8.4.2}

\[
\langle a,b,c,d\,|\,[b,a],\,[d,c]\rangle 
\]

The number of conjugacy classes is $p^{5}+3p^{4}-2p^{3}-2p^{2}+p$, and the
automorphism group has order $2(p^{2}-1)^{2}(p^{2}-p)^{2}p^{16}$.

The action of the automorphism group on $G/G^{\prime }$ is given by matrices
in GL$(4,p)$ of the form%
\[
\left( 
\begin{array}{cccc}
\ast & \ast & 0 & 0 \\ 
\ast & \ast & 0 & 0 \\ 
0 & 0 & \ast & \ast \\ 
0 & 0 & \ast & \ast%
\end{array}%
\right) \mathrm{ and }\left( 
\begin{array}{cccc}
0 & 0 & \ast & \ast \\ 
0 & 0 & \ast & \ast \\ 
\ast & \ast & 0 & 0 \\ 
\ast & \ast & 0 & 0%
\end{array}%
\right) . 
\]

\bigskip \noindent\textbf{Group 8.4.3}

\[
\langle a,b,c,d\,|\,[b,a],\,[d,b][c,a]\rangle 
\]

The number of conjugacy classes is $p^{5}+2p^{4}-p^{3}-p^{2}$, and the
automorphism group has order $(p-1)(p^{2}-1)(p^{2}-p)p^{20}$.

The action of the automorphism group on $G/G^{\prime }$ is given by matrices
in GL$(4,p)$ of the form%
\[
\left( 
\begin{array}{cccc}
\alpha & \beta & 0 & 0 \\ 
\gamma & \delta & 0 & 0 \\ 
\ast & \ast & -\lambda \delta & \lambda \gamma \\ 
\ast & \ast & \lambda \beta & -\lambda \alpha%
\end{array}%
\right) , 
\]%
with $\lambda (\alpha \delta -\beta \gamma )\neq 0$.

\bigskip \noindent\textbf{Group 8.4.4}

\[
\langle a,b,c,d\,|\,[d,b][c,a],\,[d,c][b,a]^{\omega }\rangle 
\]

The number of conjugacy classes is $p^{5}+p^{4}-p$, and the automorphism
group has order $2(p^{4}-1)(p^{4}-p^{2})p^{16}$.

The action of the automorphism group on $G/G^{\prime }$ is given by matrices
in GL$(4,p)$ of the form%
\[
\left( 
\begin{array}{cccc}
\alpha & \beta & \gamma & \delta \\ 
\varepsilon & \zeta & \eta & \theta \\ 
-\omega \theta & \omega \eta & \zeta & -\varepsilon \\ 
\omega \delta & -\omega \gamma & -\beta & \alpha%
\end{array}%
\right) \mathrm{ and }\left( 
\begin{array}{cccc}
\alpha & \beta & \gamma & \delta \\ 
\varepsilon & \zeta & \eta & \theta \\ 
\omega \theta & -\omega \eta & -\zeta & \varepsilon \\ 
-\omega \delta & \omega \gamma & \beta & -\alpha%
\end{array}%
\right) . 
\]%
In these two matrices the first row can be anything non-zero. The first row
determines the fourth row up to sign, and the second row can be anything not
in the span of rows one and four. The subgroup of matrices with first row $%
(1,0,0,0)$ is generated by matrices with second row $(\varepsilon ,\zeta
,\eta ,\theta )$ with one (or both) of $\zeta ,\eta $ non-zero. The full
matrix group is then generated by this subgroup and matrices with a general
first row $(\alpha ,\beta ,\gamma ,\delta )$, and second row $(0,1,0,0)$ if $%
\alpha \neq 0$ or $\delta \neq 0$, or second row $(1,0,0,0)$ if $\alpha
=\delta =0$.

\subsection{Five generator groups}

\noindent\textbf{Group 8.5.1}

\[
\langle
a,b,c,d,e\,|\,[e,a],\,[c,b],\,[d,b],\,[e,b],\,[d,c],\,[e,c],\,[e,d]\rangle 
\]

The number of conjugacy classes is $p^{6}+p^{5}-p^{3}$, and the automorphism
group has order $(p-1)^{2}(p^{3}-1)(p^{3}-p)(p^{3}-p^{2})p^{22}$.

The action of the automorphism group on $G/G^{\prime }$ is given by matrices
in GL$(5,p)$ of the form%
\[
\left( 
\begin{array}{ccccc}
\ast & \ast & \ast & \ast & \ast \\ 
0 & \ast & \ast & \ast & \ast \\ 
0 & \ast & \ast & \ast & \ast \\ 
0 & \ast & \ast & \ast & \ast \\ 
0 & 0 & 0 & 0 & \ast%
\end{array}%
\right) . 
\]

\bigskip \noindent\textbf{Group 8.5.2}

\[
\langle
a,b,c,d,e\,|\,[d,a],\,[e,a],\,[d,b],\,[e,b],\,[d,c],\,[e,c],\,[e,d]\rangle 
\]

The number of conjugacy classes is $p^{6}+p^{5}-p^{3}$, and the automorphism
group has order $(p^{2}-1)(p^{2}-p)(p^{3}-1)(p^{3}-p)(p^{3}-p^{2})p^{21}$.

The action of the automorphism group on $G/G^{\prime }$ is given by matrices
in GL$(5,p)$ of the form%
\[
\left( 
\begin{array}{ccccc}
\ast & \ast & \ast & \ast & \ast \\ 
\ast & \ast & \ast & \ast & \ast \\ 
\ast & \ast & \ast & \ast & \ast \\ 
0 & 0 & 0 & \ast & \ast \\ 
0 & 0 & 0 & \ast & \ast%
\end{array}%
\right) . 
\]

\bigskip \noindent\textbf{Group 8.5.3}

\[
\langle
a,b,c,d,e\,|\,[d,a],\,[e,a],\,[c,b],\,[e,b],\,[d,c],\,[e,c],\,[e,d]\rangle 
\]

The number of conjugacy classes is $3p^{5}-3p^{3}+p^{2}$, and the
automorphism group has order $2(p-1)^{5}p^{23}$.

The action of the automorphism group on $G/G^{\prime }$ is given by matrices
in GL$(5,p)$ of the form%
\[
\left( 
\begin{array}{ccccc}
\ast & 0 & \ast & \ast & \ast \\ 
0 & \ast & \ast & \ast & \ast \\ 
0 & 0 & \ast & 0 & \ast \\ 
0 & 0 & 0 & \ast & \ast \\ 
0 & 0 & 0 & 0 & \ast%
\end{array}%
\right) \mathrm{ and }\left( 
\begin{array}{ccccc}
0 & \ast & \ast & \ast & \ast \\ 
\ast & 0 & \ast & \ast & \ast \\ 
0 & 0 & 0 & \ast & \ast \\ 
0 & 0 & \ast & 0 & \ast \\ 
0 & 0 & 0 & 0 & \ast%
\end{array}%
\right) . 
\]

\bigskip \noindent\textbf{Group 8.5.4}

\[
\langle
a,b,c,d,e\,|\,[d,a]=[c,b],\,[e,a],\,[d,b],\,[e,b],\,[d,c],\,[e,c],\,[e,d]%
\rangle 
\]

The number of conjugacy classes is $2p^{5}+p^{4}-2p^{3}$, and the
automorphism group has order $(p-1)^{2}(p^{2}-1)(p^{2}-p)p^{24}$.

The action of the automorphism group on $G/G^{\prime }$ is given by matrices
in GL$(5,p)$ of the form%
\[
\left( 
\begin{array}{ccccc}
\lambda & \ast & \ast & \ast & \ast \\ 
0 & \alpha & \beta & \ast & \ast \\ 
0 & \gamma & \delta & \ast & \ast \\ 
0 & 0 & 0 & \lambda ^{-1}(\alpha \delta -\beta \gamma ) & \ast \\ 
0 & 0 & 0 & 0 & \ast%
\end{array}%
\right) , 
\]%
with $\lambda (\alpha \delta -\beta \gamma )\neq 0$.

\bigskip \noindent\textbf{Group 8.5.5}

\[
\langle
a,b,c,d,e\,|\,[c,a],\,[d,a]=[c,b],\,[e,a],\,[d,b],\,[e,b],\,[e,c],\,[e,d]%
\rangle 
\]

The number of conjugacy classes is $2p^{5}+p^{4}-2p^{3}$, and the
automorphism group has order $(p^{2}-1)^{2}(p^{2}-p)^{2}p^{19}$.

The action of the automorphism group on $G/G^{\prime }$ is given by matrices
in GL$(5,p)$ of the form%
\[
\left( 
\begin{array}{ccccc}
\lambda \alpha & \lambda \beta & \mu \alpha & -\mu \beta & \ast \\ 
-\lambda \gamma & -\lambda \delta & -\mu \gamma & \mu \delta & \ast \\ 
\nu \alpha & \nu \beta & \xi \alpha & -\xi \beta & \ast \\ 
\nu \gamma & \nu \delta & \xi \gamma & -\xi \delta & \ast \\ 
0 & 0 & 0 & 0 & \ast%
\end{array}%
\right) , 
\]%
with $(\alpha \delta -\beta \gamma )(\lambda \xi -\mu \nu )\neq 0$

\bigskip \noindent\textbf{Group 8.5.6}

\[
\langle a,b,c,d,e\,|\,[d,b]=[c,a]^{\omega
},\,[d,c]=[b,a],\,[e,a],\,[c,b],\,[e,b],\,[e,c],\,[e,d]\rangle 
\]

The number of conjugacy classes is $p^{5}+2p^{4}-p^{3}-p^{2}$, and the
automorphism group has order $2(p-1)(p^{2}-1)^{2}p^{23}$.

The action of the automorphism group on $G/G^{\prime }$ is given by matrices
in GL$(5,p)$ of the form%
\[
\left( 
\begin{array}{ccccc}
\alpha & \ast & \ast & \beta & \ast \\ 
0 & \gamma & \omega \delta & 0 & \ast \\ 
0 & \delta & \gamma & 0 & \ast \\ 
\omega \beta & \ast & \ast & \alpha & \ast \\ 
0 & 0 & 0 & 0 & \ast%
\end{array}%
\right) \mathrm{ and }\left( 
\begin{array}{ccccc}
\alpha & \ast & \ast & \beta & \ast \\ 
0 & \gamma & \omega \delta & 0 & \ast \\ 
0 & -\delta & -\gamma & 0 & \ast \\ 
-\omega \beta & \ast & \ast & -\alpha & \ast \\ 
0 & 0 & 0 & 0 & \ast%
\end{array}%
\right) , 
\]%
with $\alpha $ and $\beta $ not both zero, and with $\gamma $ and $\delta $
not both zero.

\bigskip \noindent\textbf{Group 8.5.7}

\[
\langle
a,b,c,d,e\,|\,[e,b]=[c,a][d,b]^{m},\,[d,c]=[b,a],\,[e,c]=[d,b],\,[d,a],%
\,[e,a],\,[c,b],\,[e,d]\rangle , 
\]%
where $m$ is chosen so that $x^{3}+mx-1$ is irreducible over GF$(p)$.
Different choices of $m$ give isomorphic groups. Note that the discriminant
of $x^{3}+mx-1$ is $-4m^{3}-27$ and this must be a square in GF$(p)$ --- we
let $u^{2}=-4m^{3}-27$.

The number of conjugacy classes is $p^{5}+p^{4}-p$, and the automorphism
group has order $3(p-1)(p^{3}-1)p^{18}$.

The action of the automorphism group on $G/G^{\prime }$ is given by matrices
in GL$(5,p)$ of the form%
\[
\left( 
\begin{array}{ccccc}
\ast & 0 & 0 & \ast & \ast \\ 
\ast & \ast & \ast & \ast & \ast \\ 
\ast & \ast & \ast & \ast & \ast \\ 
\ast & 0 & 0 & \ast & \ast \\ 
\ast & 0 & 0 & \ast & \ast%
\end{array}%
\right) 
\]%
where the first row can be anything non-zero of the form shown ($(p^{3}-1)$
possibilities). If we restrict ourselves to matrices with first row $%
(1,0,0,0,0)$ then we get a subgroup of GL$(5,p)$ of order $3(p-1)p^{3}$
generated by matrices of the following form%
\[
\left( 
\begin{array}{ccccc}
1 & 0 & 0 & 0 & 0 \\ 
\alpha & \beta & 0 & \gamma & \delta \\ 
-\gamma & 0 & \beta & m\alpha +\delta & -\alpha \\ 
0 & 0 & 0 & 1 & 0 \\ 
0 & 0 & 0 & 0 & 1%
\end{array}%
\right) 
\]%
with $\beta \neq 0$, and%
\[
\left( 
\begin{array}{ccccc}
1 & 0 & 0 & 0 & 0 \\ 
0 & \frac{u-9}{2m^{2}} & \frac{3}{m} & 0 & 0 \\ 
0 & 1 & \frac{u+9}{2m^{2}} & 0 & 0 \\ 
-\frac{(u+3)m}{2u} & 0 & 0 & \frac{u^{2}-2u+9}{4u} & \frac{m^{2}}{u} \\ 
-\frac{(u+1)m^{2}}{2u} & 0 & 0 & \frac{(u^{2}+3)m}{4u} & -\frac{u^{2}+2u+9}{%
4u}%
\end{array}%
\right) . 
\]%
(The cube of the second matrix has the form of the first, with $\alpha
=\gamma =\delta =0$, and $\beta =4u^{3}m^{-3}$.) The full subgroup of GL$%
(5,p)$ giving the action of the automorphism group on $G/G^{\prime }$ is
generated by the matrices above together with $p^{3}-1$ matrices of the form%
\[
\left( 
\begin{array}{ccccc}
\alpha & 0 & 0 & \beta & \gamma \\ 
0 & 1 & 0 & 0 & 0 \\ 
0 & 0 & 1 & 0 & 0 \\ 
\gamma & 0 & 0 & \alpha & -\beta -m\gamma \\ 
\beta +m\gamma & 0 & 0 & -\gamma & \alpha -m\beta -m^{2}\gamma%
\end{array}%
\right) , 
\]%
where $\alpha ,\beta ,\gamma $ are not all zero.

\newpage \noindent\textbf{Group 8.5.8}

\[
\langle a,b,c,d,e\,|\,[d,a][c,b]^{\omega
},\,[e,a],\,[e,b]=[c,a],\,[d,b],\,[d,c]=[b,a],\,[e,c],\,[e,d]\rangle 
\]

The number of conjugacy classes is $p^{5}+p^{4}-p$, and the automorphism
group has order $2(p-1)(p^{2}-1)p^{18}$.

The action of the automorphism group on $G/G^{\prime }$ is given by a group $%
H$ of matrices in GL$(5,p)$ of the form%
\[
\left( 
\begin{array}{ccccc}
\ast & 0 & \ast & 0 & \ast \\ 
\ast & \ast & \ast & \ast & \ast \\ 
\ast & 0 & \ast & 0 & \ast \\ 
\ast & 0 & \ast & \ast & \ast \\ 
0 & 0 & 0 & 0 & \ast%
\end{array}%
\right) . 
\]%
The group $H$ has a subgroup of order $p^{3}(p-1)$ consisting of matrices of
the form%
\[
\left( 
\begin{array}{ccccc}
1 & 0 & 0 & 0 & 2\mu \xi \\ 
\lambda & \xi ^{-1} & \mu & 0 & \nu \\ 
0 & 0 & 1 & 0 & -2\lambda \xi \\ 
\mu & 0 & \omega \lambda & \xi ^{-1} & -\omega \lambda ^{2}\xi +\mu ^{2}\xi
\\ 
0 & 0 & 0 & 0 & \xi%
\end{array}%
\right) 
\]%
with $\xi \neq 0$.

If we premultiply a general matrix%
\[
\left( 
\begin{array}{ccccc}
\ast & 0 & \ast & 0 & \ast \\ 
\ast & \ast & \ast & \ast & \ast \\ 
\ast & 0 & \ast & 0 & \ast \\ 
\ast & 0 & \ast & \ast & \ast \\ 
0 & 0 & 0 & 0 & \ast%
\end{array}%
\right) 
\]%
in $H$ by a suitable matrix from this subgroup of order $p^{3}(p-1)$ we can
obtain a matrix of the form%
\[
\left( 
\begin{array}{ccccc}
\ast & 0 & \ast & 0 & \ast \\ 
0 & \ast & 0 & \ast & 0 \\ 
\ast & 0 & \ast & 0 & \ast \\ 
\ast & 0 & \ast & \ast & \ast \\ 
0 & 0 & 0 & 0 & \ast%
\end{array}%
\right) , 
\]%
and the most general matrices of this form arising in the action of the
automorphism group on $G/G^{\prime }$ are

\[
\xi \left( 
\begin{array}{ccccc}
\alpha & 0 & \omega \beta & 0 & 0 \\ 
0 & \alpha ^{2}-\omega \beta ^{2} & 0 & 0 & 0 \\ 
\beta & 0 & \alpha & 0 & 0 \\ 
0 & 0 & 0 & \alpha ^{2}-\omega \beta ^{2} & 0 \\ 
0 & 0 & 0 & 0 & 1%
\end{array}%
\right) 
\]%
and 
\[
\xi \left( 
\begin{array}{ccccc}
\alpha & 0 & \omega \beta & 0 & 0 \\ 
0 & -\alpha ^{2}+\omega \beta ^{2} & 0 & 0 & 0 \\ 
-\beta & 0 & -\alpha & 0 & 0 \\ 
0 & 0 & 0 & \alpha ^{2}-\omega \beta ^{2} & 0 \\ 
0 & 0 & 0 & 0 & 1%
\end{array}%
\right) . 
\]%
There are $2(p-1)(p^{2}-1)$ of these matrices, and so the order of the
subgroup of GL$(5,p)$ generated by all these matrices is $%
2p^{3}(p-1)(p^{2}-1)$.

\bigskip \noindent\textbf{Group 8.5.9}

\[
\langle
a,b,c,d,e\,|\,[d,a],\,[e,a],\,[e,b]=[c,a],\,[d,b],\,[d,c]=[b,a],\,[e,c],%
\,[e,d]=[c,b]\rangle 
\]

The number of conjugacy classes is $p^{5}+p^{4}-p$, and the automorphism
group has order $(p+1)(p-1)^{2}p^{16}$.

The action of the automorphism group on $G/G^{\prime }$ is given by a
subgroup $H$ of GL$(5,p)$ of order $(p+1)(p-1)^{2}p$. There is a subgroup $K$
of $H$ of order $p(p-1)^{2}$ consisting of matrices of the form%
\[
\lambda \left( 
\begin{array}{ccccc}
\alpha ^{2} & \frac{1}{2}\alpha \beta & 0 & \frac{1}{4}\beta ^{2} & 0 \\ 
0 & \alpha & 0 & \beta & 0 \\ 
-\frac{3}{2}\alpha ^{2}\beta & -\frac{3}{8}\alpha \beta ^{2} & \alpha ^{3} & 
-\frac{1}{8}\beta ^{3} & 0 \\ 
0 & 0 & 0 & 1 & 0 \\ 
\frac{3}{8}\alpha ^{2}\beta ^{2} & \frac{1}{16}\alpha \beta ^{3} & -\frac{1}{%
2}\alpha ^{3}\beta & \frac{1}{64}\beta ^{4} & \alpha ^{4}%
\end{array}%
\right) 
\]
$\allowbreak $with $\alpha ,\lambda \neq 0$, and $K$ consists of all the
matrices in $H$ with fourth row a scalar multiple of $(0,0,0,1,0)$. When $%
p=3 $ there is a matrix%
\[
A=\left( 
\begin{array}{ccccc}
1 & 0 & 1 & 0 & 1 \\ 
0 & 1 & 0 & 0 & 2 \\ 
0 & 0 & 1 & 0 & 2 \\ 
0 & 2 & 2 & 1 & 1 \\ 
0 & 0 & 0 & 0 & 1%
\end{array}%
\right) 
\]%
in $H$, and when $p\neq 3$, there is a matrix%
\[
A=\left( 
\begin{array}{ccccc}
-\frac{1}{3} & 0 & 0 & \frac{3}{2} & \frac{2}{27} \\ 
0 & 1 & \frac{2}{9} & 3 & -\frac{4}{27} \\ 
0 & \frac{1}{2} & \frac{1}{9} & -\frac{3}{2} & \frac{2}{27} \\ 
1 & 1 & -\frac{2}{9} & \frac{3}{2} & \frac{2}{27} \\ 
\frac{1}{4} & -\frac{1}{4} & \frac{1}{18} & \frac{3}{8} & \frac{1}{54}%
\end{array}%
\right) 
\]%
in $H$. In both cases $H$ is generated by $A$ and $K$. (There are $(p^{2}-1)$
possibilities for row 4, and all these possibilities arise as the fourth row
in elements of $K$ or $AK$.)

\bigskip \noindent\textbf{Group 8.5.10}

\[
\langle
a,b,c,d,e\,|\,[d,a],\,[e,a],\,[d,b],\,[e,b]=[c,a],\,[d,c]=[b,a],\,[e,c],%
\,[e,d]\rangle 
\]

The number of conjugacy classes is $p^{5}+2p^{4}-p^{3}-p^{2}$, and the
automorphism group has order $(p-1)(p^{2}-1)(p^{2}-p)p^{21}$.

The action of the automorphism group on $G/G^{\prime }$ is given by matrices
in GL$(5,p)$ of the form%
\[
\left( 
\begin{array}{ccccc}
\ast & 0 & 0 & \ast & \ast \\ 
\ast & \alpha & \beta & \ast & \ast \\ 
\ast & \gamma & \delta & \ast & \ast \\ 
\ast & 0 & 0 & \ast & \ast \\ 
\ast & 0 & 0 & \ast & \ast%
\end{array}%
\right) 
\]%
with $\alpha \delta -\beta \gamma \neq 0$. Rows 2 and 3 are arbitrary
subject to the condition that $\alpha \delta -\beta \gamma \neq 0$, and rows
2 and 3 then determine rows 1, 4 and 5 up to multiplication by a scalar.
(The same scalar for each of the three rows. There is a subgroup of this
group of matrices of order $(p-1)p^{6}$ consisting of matrices of the form%
\[
\left( 
\begin{array}{ccccc}
\lambda & 0 & 0 & 0 & 0 \\ 
\ast & 1 & 0 & \ast & \ast \\ 
\ast & 0 & 1 & \ast & \ast \\ 
0 & 0 & 0 & \lambda & 0 \\ 
0 & 0 & 0 & 0 & \lambda%
\end{array}%
\right) \;(\lambda \neq 0), 
\]%
and the full subgroup of GL$(5,p)$ giving the action of the automorphism
group on $G/G^{\prime }$ is generated by these matrices, together with the
matrices%
\[
\left( 
\begin{array}{ccccc}
\lambda & 0 & 0 & 0 & 0 \\ 
0 & \lambda & 0 & 0 & 0 \\ 
0 & 0 & 1 & 0 & 0 \\ 
0 & 0 & 0 & \lambda ^{2} & 0 \\ 
0 & 0 & 0 & 0 & 1%
\end{array}%
\right) \mathrm{ }(\lambda \neq 0)\mathrm{ and }\left( 
\begin{array}{ccccc}
1 & 0 & 0 & -2 & 0 \\ 
0 & -1 & 1 & 0 & 0 \\ 
0 & -1 & 0 & 0 & 0 \\ 
1 & 0 & 0 & -1 & -1 \\ 
0 & 0 & 0 & -1 & 0%
\end{array}%
\right) . 
\]

\bigskip \noindent\textbf{Group 8.5.11}

\[
\langle
a,b,c,d,e\,|\,[d,a],\,[e,a],\,[c,b],\,[e,b]=[c,a],\,[d,c]=[b,a],\,[e,c],%
\,[e,d]\rangle 
\]

The number of conjugacy classes is $p^{5}+2p^{4}-p^{3}-p^{2}$, and the
automorphism group has order $(p-1)^{3}p^{21}$.

The action of the automorphism group on $G/G^{\prime }$ is given by matrices
in GL$(5,p)$ of the form%
\[
\lambda \left( 
\begin{array}{ccccc}
\varepsilon \eta & 0 & 0 & 0 & \beta \varepsilon \eta +\varepsilon \zeta \\ 
\alpha & 1 & \beta & \gamma & \delta \\ 
-\gamma \varepsilon & 0 & \varepsilon & 0 & -\alpha \varepsilon -\beta
\gamma \varepsilon \\ 
\zeta & 0 & 0 & \eta & \theta \\ 
0 & 0 & 0 & 0 & \varepsilon ^{2}\eta%
\end{array}%
\right) 
\]%
with $\lambda ,\varepsilon ,\eta \neq 0$.

\bigskip \noindent\textbf{Group 8.5.12}

\[
\langle
a,b,c,d,e\,|\,[e,a],\,[c,b],\,[d,b],\,[e,c],\,[e,d],\,[d,c]=[b,a],%
\,[e,b]=[c,a]\rangle 
\]

The number of conjugacy classes is $p^{5}+2p^{4}-p^{3}-p^{2}$, and the
automorphism group has order $(p-1)^{3}p^{19}$.

The action of the automorphism group on $G/G^{\prime }$ is given by a group $%
H$ of matrices in GL$(5,p)$ where $H$ is generated by matrices of the form%
\[
\lambda \left( 
\begin{array}{ccccc}
\alpha & 0 & 0 & 0 & 0 \\ 
0 & \beta & 0 & 0 & 0 \\ 
0 & 0 & \alpha ^{-1}\beta & 0 & 0 \\ 
0 & 0 & 0 & \alpha ^{2} & 0 \\ 
0 & 0 & 0 & 0 & 1%
\end{array}%
\right) \mathrm{ }(\alpha ,\beta ,\lambda \neq 0)\mathrm{ and }\left( 
\begin{array}{ccccc}
1 & 0 & \ast & 0 & 2\gamma \\ 
0 & 1 & \gamma & 0 & \ast \\ 
0 & 0 & 1 & 0 & 0 \\ 
\gamma & 0 & \ast & 1 & \gamma ^{2} \\ 
0 & 0 & 0 & 0 & 1%
\end{array}%
\right) . 
\]

\bigskip \noindent\textbf{Group 8.5.13}

\[
\langle a,b,c,d,e\,|\,[d,a],\,[e,a],\,[c,b],\,[e,b][d,b]^{\omega
}=[c,a],\,[d,c]=[b,a],\,[e,c],\,[e,d]\rangle 
\]

The number of conjugacy classes is $p^{5}+2p^{4}-p^{3}-p^{2}$, and the
automorphism group has order $2(p-1)^{2}(p^{2}-1)p^{18}$.

The action of the automorphism group on $G/G^{\prime }$ is given by matrices
in GL$(5,p)$ of the form%
\[
\lambda \left( 
\begin{array}{ccccc}
\alpha & 0 & 0 & \omega \beta & \beta \\ 
-\varepsilon & 1 & 0 & -\delta & \gamma \\ 
\delta & 0 & 1 & \omega \varepsilon & \varepsilon \\ 
\beta & 0 & 0 & \alpha & \zeta \\ 
0 & 0 & 0 & 0 & \alpha -\omega \zeta%
\end{array}%
\right) \mathrm{ and }\lambda \left( 
\begin{array}{ccccc}
\alpha & 0 & 0 & \omega \beta & \beta \\ 
\varepsilon & -1 & 0 & \delta & -\gamma \\ 
\delta & 0 & 1 & \omega \varepsilon & \varepsilon \\ 
-\beta & 0 & 0 & -\alpha & -\zeta \\ 
0 & 0 & 0 & 0 & -\alpha +\omega \zeta%
\end{array}%
\right) , 
\]%
where $\alpha $ and $\beta $ are not both zero and $\zeta \neq \alpha \omega
^{-1}$, $\lambda \neq 0$.

\bigskip \noindent\textbf{Group 8.5.14}

\[
\langle
a,b,c,d,e\,|\,[c,b],\,[d,b],\,[e,c],\,[e,d],\,[d,c]=[b,a],\,[e,b]=[c,a],%
\,[e,a]=[d,a]^{\omega }\rangle 
\]

The number of conjugacy classes is $p^{5}+2p^{4}-p^{3}-p^{2}$, and the
automorphism group has order $2(p-1)(p^{2}-1)p^{17}$.

The action of the automorphism group on $G/G^{\prime }$ is given by a
subgroup $H=KL$ of GL$(5,p)$, where $K$ is a group of matrices of order $%
2(p-1)(p^{2}-1)$ consisting of matrices%
\[
\left( 
\begin{array}{ccccc}
\lambda (\alpha ^{2}+\omega \beta ^{2}) & 0 & 0 & 2\omega \lambda \alpha
\beta & 2\lambda \alpha \beta \\ 
0 & \alpha & \beta & 0 & 0 \\ 
0 & \omega \beta & \alpha & 0 & 0 \\ 
\lambda \alpha \beta & 0 & 0 & \lambda \alpha ^{2} & \lambda \beta ^{2} \\ 
\omega \lambda \alpha \beta & 0 & 0 & \omega ^{2}\lambda \beta ^{2} & 
\lambda \alpha ^{2}%
\end{array}%
\right) 
\]
and 
\[
\left( 
\begin{array}{ccccc}
\lambda (\alpha ^{2}+\omega \beta ^{2}) & 0 & 0 & 2\omega \lambda \alpha
\beta & 2\lambda \alpha \beta \\ 
0 & \alpha & \beta & 0 & 0 \\ 
0 & -\omega \beta & -\alpha & 0 & 0 \\ 
-\lambda \alpha \beta & 0 & 0 & -\lambda \alpha ^{2} & -\lambda \beta ^{2}
\\ 
-\omega \lambda \alpha \beta & 0 & 0 & -\omega ^{2}\lambda \beta ^{2} & 
-\lambda \alpha ^{2}%
\end{array}%
\right) , 
\]%
with $\lambda \neq 0$ and $\alpha ,\beta $ not both zero, and where $L$ is a
normal subgroup of $H$ of order $p^{2}$ consisting of matrices%
\[
\left( 
\begin{array}{ccccc}
1 & \omega \alpha & \beta & 0 & 0 \\ 
0 & 1 & 0 & 0 & 0 \\ 
0 & 0 & 1 & 0 & 0 \\ 
0 & 0 & \alpha & 1 & 0 \\ 
0 & \omega \beta & 0 & 0 & 1%
\end{array}%
\right) . 
\]

\bigskip \noindent\textbf{Group 8.5.15}

\[
\langle
a,b,c,d,e\,|\,[d,a],\,[e,a],\,[c,b],\,[d,b],\,[e,c],\,[d,c]=[b,a],%
\,[e,b]=[c,a]\rangle 
\]

The number of conjugacy classes is $p^{5}+2p^{4}-p^{3}-p^{2}$, and the
automorphism group has order $2(p-1)^{3}p^{17}$.

The action of the automorphism group on $G/G^{\prime }$ is given by matrices
in GL$(5,p)$ of the form%
\[
\lambda \left( 
\begin{array}{ccccc}
1 & \alpha & \beta & 0 & 0 \\ 
0 & \gamma \delta ^{-1} & 0 & 0 & 0 \\ 
0 & 0 & \gamma & 0 & 0 \\ 
0 & -\beta \delta ^{-1} & 0 & \delta ^{-1} & 0 \\ 
0 & 0 & -\alpha \delta & 0 & \delta%
\end{array}%
\right) \mathrm{ and }\lambda \left( 
\begin{array}{ccccc}
1 & \alpha & \beta & 0 & 0 \\ 
0 & 0 & \gamma & 0 & 0 \\ 
0 & \gamma \delta ^{-1} & 0 & 0 & 0 \\ 
0 & 0 & -\beta \delta & 0 & \delta \\ 
0 & -\beta \delta ^{-1} & 0 & \delta ^{-1} & 0%
\end{array}%
\right) , 
\]%
with $\gamma ,\delta ,\lambda \neq 0$.

\bigskip \noindent\textbf{Group 8.5.16}

\[
\langle
a,b,c,d,e\,|\,[d,a],\,[e,a],\,[c,b],\,[d,b],\,[e,b]=[c,a],\,[e,c],\,[e,d]%
\rangle 
\]

The number of conjugacy classes is $p^{5}+3p^{4}-2p^{3}-2p^{2}+p$, and the
automorphism group has order $(p-1)^{4}p^{19}$.

The action of the automorphism group on $G/G^{\prime }$ is given by matrices
in GL$(5,p)$ of the form%
\[
\left( 
\begin{array}{ccccc}
\alpha & 0 & 0 & 0 & 0 \\ 
0 & \beta & 0 & 0 & 0 \\ 
0 & 0 & \alpha ^{-1}\beta \delta & 0 & 0 \\ 
0 & 0 & 0 & \gamma & 0 \\ 
0 & 0 & 0 & 0 & \delta%
\end{array}%
\right) \left( 
\begin{array}{ccccc}
1 & 0 & 0 & 0 & \varepsilon \\ 
-\theta & 1 & 0 & 0 & \zeta \\ 
0 & 0 & 1 & \eta & \theta \\ 
0 & 0 & 0 & 1 & 0 \\ 
0 & 0 & 0 & 0 & 1%
\end{array}%
\right) , 
\]%
with $\alpha ,\beta ,\gamma ,\delta \neq 0$.

\bigskip \noindent\textbf{Group 8.5.17}

\[
\langle
a,b,c,d,e\,|\,[e,b]=[c,a],\,[d,c]=[b,a][d,a]^{-1},\,[e,a],\,[c,b],\,[d,b],%
\,[e,c],\,[e,d]\rangle 
\]

The number of conjugacy classes is $p^{5}+3p^{4}-2p^{3}-2p^{2}+p$, and the
automorphism group has order $2(p-1)^{3}p^{18}$.

The action of the automorphism group on $G/G^{\prime }$ is given by matrices
in GL$(5,p)$ of the form%
\[
\lambda \left( 
\begin{array}{ccccc}
\alpha & 0 & 0 & 0 & 0 \\ 
0 & \alpha ^{2} & 0 & 0 & 0 \\ 
0 & 0 & \alpha & 0 & 0 \\ 
0 & 0 & 0 & \alpha ^{2} & 0 \\ 
0 & 0 & 0 & 0 & 1%
\end{array}%
\right) \left( 
\begin{array}{ccccc}
1 & 0 & 1-\gamma & 0 & \varepsilon \\ 
\beta & \gamma & \beta & 0 & \delta \\ 
0 & 0 & \gamma & 0 & -\beta \\ 
\beta & 0 & \beta +\beta \gamma -\gamma \varepsilon & \gamma & \beta
\varepsilon -\beta ^{2} \\ 
0 & 0 & 0 & 0 & 1%
\end{array}%
\right) 
\]%
and%
\[
\lambda \left( 
\begin{array}{ccccc}
\alpha & 0 & 0 & 0 & 0 \\ 
0 & \alpha ^{2} & 0 & 0 & 0 \\ 
0 & 0 & \alpha & 0 & 0 \\ 
0 & 0 & 0 & \alpha ^{2} & 0 \\ 
0 & 0 & 0 & 0 & 1%
\end{array}%
\right) \left( 
\begin{array}{ccccc}
1 & 0 & 1+\gamma & 0 & \varepsilon \\ 
0 & \gamma & -\beta \gamma & -\gamma & \delta \\ 
-1 & 0 & -1 & 0 & \beta \\ 
\beta +\varepsilon & 0 & \beta +\varepsilon -\beta \gamma & -\gamma & -\beta
\varepsilon -\beta ^{2} \\ 
0 & 0 & 0 & 0 & 1%
\end{array}%
\right) , 
\]%
with $\alpha ,\gamma ,\lambda \neq 0$.

\bigskip \noindent\textbf{$\allowbreak $Group 8.5.18}

\[
\langle
a,b,c,d,e\,|\,[e,c][e,b],\,[c,a],\,[d,a],\,[e,a],\,[c,b],\,[d,b],\,[e,d]%
\rangle 
\]

The number of conjugacy classes is $p^{5}+4p^{4}-3p^{3}-3p^{2}+2p$, and the
automorphism group has order $6(p-1)^{4}p^{18}$.

The action of the automorphism group on $G/G^{\prime }$ is generate by
matrices in GL$(5,p)$ of the form%
\[
\left( 
\begin{array}{ccccc}
\alpha & 0 & 0 & 0 & 0 \\ 
0 & \beta & 0 & 0 & 0 \\ 
0 & 0 & \beta & 0 & 0 \\ 
0 & 0 & 0 & \gamma & 0 \\ 
0 & 0 & 0 & 0 & \delta%
\end{array}%
\right) \left( 
\begin{array}{ccccc}
1 & 0 & 0 & 0 & 0 \\ 
\varepsilon & 1 & 0 & 0 & \zeta \\ 
0 & 0 & 1 & \eta & -\zeta \\ 
0 & 0 & 0 & 1 & 0 \\ 
0 & 0 & 0 & 0 & 1%
\end{array}%
\right) \;(\alpha ,\beta ,\gamma ,\delta \neq 0), 
\]%
\[
\left( 
\begin{array}{ccccc}
1 & 0 & 0 & 0 & 0 \\ 
0 & 1 & 1 & 0 & 0 \\ 
0 & 0 & -1 & 0 & 0 \\ 
0 & 0 & 0 & 0 & 1 \\ 
0 & 0 & 0 & 1 & 0%
\end{array}%
\right) \mathrm{ and }\left( 
\begin{array}{ccccc}
0 & 0 & 0 & 1 & 0 \\ 
0 & 0 & 1 & 0 & 0 \\ 
0 & 1 & 0 & 0 & 0 \\ 
1 & 0 & 0 & 0 & 0 \\ 
0 & 0 & 0 & 0 & 1%
\end{array}%
\right) . 
\]%
Here the matrices of the first kind above form a normal subgroup of order $%
(p-1)^{4}p^{3}$.

\bigskip \noindent\textbf{Group 8.5.19}

\[
\langle
a,b,c,d,e\,|\,[d,a],\,[e,a],\,[c,b],\,[d,b],\,[e,b],\,[d,c],\,[e,c]\rangle 
\]

The number of conjugacy classes is $2p^{5}+2p^{4}-3p^{3}-p^{2}+p$, and the
automorphism group has order $(p-1)(p^{2}-1)^{2}(p^{2}-p)^{2}p^{17}$.

The action of the automorphism group on $G/G^{\prime }$ is given by matrices
in GL$(5,p)$ of the form%
\[
\left( 
\begin{array}{ccccc}
\ast & \ast & \ast & 0 & 0 \\ 
0 & \ast & \ast & 0 & 0 \\ 
0 & \ast & \ast & 0 & 0 \\ 
0 & 0 & 0 & \ast & \ast \\ 
0 & 0 & 0 & \ast & \ast%
\end{array}%
\right) . 
\]

\bigskip \noindent\textbf{Group 8.5.20}

\[
\langle
a,b,c,d,e\,|\,[d,c]=[b,a],\,[c,a],\,[d,a],\,[e,a],\,[d,b],\,[e,b],\,[e,d]%
\rangle 
\]

The number of conjugacy classes is $2p^{5}+p^{4}-2p^{3}$, and the
automorphism group has order $(p-1)^{4}p^{23}$.

The action of the automorphism group on $G/G^{\prime }$ is given by matrices
in GL$(5,p)$ of the form%
\[
\left( 
\begin{array}{ccccc}
\alpha  & 0 & 0 & 0 & 0 \\ 
0 & \beta  & 0 & 0 & 0 \\ 
0 & 0 & \gamma  & 0 & 0 \\ 
0 & 0 & 0 & \alpha \beta \gamma ^{-1} & 0 \\ 
0 & 0 & 0 & 0 & \delta 
\end{array}%
\right) \left( 
\begin{array}{ccccc}
1 & 0 & 0 & \varepsilon  & 0 \\ 
\ast  & 1 & 0 & \ast  & \ast  \\ 
\ast  & \varepsilon  & 1 & \ast  & \ast  \\ 
0 & 0 & 0 & 1 & 0 \\ 
0 & 0 & 0 & \ast  & 1%
\end{array}%
\right) ,
\]%
with $\alpha ,\beta ,\gamma ,\delta \neq 0$.\newpage

\bigskip \noindent\textbf{Group 8.5.21}

\[
\langle
a,b,c,d,e\,|\,[d,c]=[b,a],\,[c,a],\,[d,a],\,[e,a],\,[c,b],\,[d,b],\,[e,d]%
\rangle 
\]

The number of conjugacy classes is $2p^{5}+p^{4}-2p^{3}$, and the
automorphism group has order $(p-1)^{2}(p^{2}-1)(p^{2}-p)p^{20}$.

The action of the automorphism group on $G/G^{\prime }$ is given by the
subgroup of GL$(5,p)$ generated by matrices of the form%
\[
\left( 
\begin{array}{ccccc}
\alpha & 0 & 0 & 0 & 0 \\ 
0 & \beta & 0 & 0 & 0 \\ 
0 & 0 & \gamma & 0 & 0 \\ 
0 & 0 & 0 & \alpha \beta \gamma ^{-1} & 0 \\ 
0 & 0 & 0 & 0 & \delta%
\end{array}%
\right) \;(\alpha ,\beta ,\gamma ,\delta \neq 0) 
\]%
and%
\[
\left( 
\begin{array}{ccccc}
1 & 0 & 0 & 0 & 0 \\ 
\ast & 1 & 0 & \varepsilon & 0 \\ 
-\varepsilon & 0 & 1 & \ast & 0 \\ 
0 & 0 & 0 & 1 & 0 \\ 
\ast & 0 & 0 & \ast & 1%
\end{array}%
\right) ,\;\left( 
\begin{array}{ccccc}
-1 & 0 & 0 & 1 & 0 \\ 
\varepsilon & 0 & -1 & \zeta & 0 \\ 
\ast & 1 & -1 & \varepsilon +\zeta & 0 \\ 
-1 & 0 & 0 & 0 & 0 \\ 
\ast & 0 & 0 & \ast & 1%
\end{array}%
\right) 
\]%
(If we let $\lambda ,\mu ,\nu ,\xi $ be the $(1,1)$, $(1,4)$, $(4,1)$, $%
(4,4) $ entries in a general matrix in this group of matrices, then $\left( 
\begin{array}{cc}
\lambda & \mu \\ 
\nu & \xi%
\end{array}%
\right) $ takes all values in GL$(2,p)$.)

\bigskip \noindent\textbf{Group 8.5.22}

\[
\langle
a,b,c,d,e\,|\,[d,c]=[b,a],\,[c,a],\,[d,a],\,[e,a],\,[c,b],\,[d,b],\,[e,b]%
\rangle 
\]

The number of conjugacy classes is $2p^{5}+p^{4}-2p^{3}$, and the
automorphism group has order $(p^{2}-1)^{2}(p^{2}-p)^{2}p^{17}$.

The action of the automorphism group on $G/G^{\prime }$ is given by the
subgroup of GL$(5,p)$ consisting of matrices of the form%
\[
\left( 
\begin{array}{ccccc}
\ast & \ast & 0 & 0 & 0 \\ 
\ast & \ast & 0 & 0 & 0 \\ 
0 & 0 & \ast & \ast & 0 \\ 
0 & 0 & \ast & \ast & 0 \\ 
0 & 0 & \ast & \ast & \ast%
\end{array}%
\right) , 
\]%
with the restriction that if we let $A$ be the elements of GL$(2,p)$ in
positions $(1,1)$, $(1,2)$, $(2,1)$, $(2,2)$ and if we let $B$ be the
element of GL$(2,p)$ in positions $(3,3)$, $(3,4)$, $(4,3)$, $(4,4)$, then $%
\det A=\det B$.

\subsection{Six generator groups}

For all these groups we take the generators to be $a,b,c,d,e,f$, and we just
give the relations, with the class two and exponent $p$ conditions
understood.

\bigskip \noindent\textbf{Group 8.6.1}

\begin{eqnarray*}
&&\lbrack c,b],\,[d,a],\,[d,b],\,[d,c],\,[e,a],\,[e,b],\,[e,c], \\
&&\lbrack e,d],\,[f,a],\,[f,b],\,[f,c],\,[f,d],\,[f,e]
\end{eqnarray*}

The number of conjugacy classes is $2p^{6}-p^{4}$ and the order of the
automorphism group is $%
(p-1)(p^{2}-1)(p^{2}-p)(p^{3}-1)(p^{3}-p)(p^{3}-p^{2})p^{23}.$

The action of the automorphism group on $G/G^{\prime }$ is given by matrices
in GL$(6,p)$ of the form%
\[
\left( 
\begin{array}{cccccc}
\ast & \ast & \ast & \ast & \ast & \ast \\ 
0 & \ast & \ast & \ast & \ast & \ast \\ 
0 & \ast & \ast & \ast & \ast & \ast \\ 
0 & 0 & 0 & \ast & \ast & \ast \\ 
0 & 0 & 0 & \ast & \ast & \ast \\ 
0 & 0 & 0 & \ast & \ast & \ast%
\end{array}%
\right) . 
\]

\bigskip \noindent\textbf{Group 8.6.2}

\begin{eqnarray*}
&&[c,b],\,[d,a],\,[d,b]=[b,a],\,[d,c],\,[e,a],\,[e,b],\,[e,c], \\
&&[e,d],\,[f,a],\,[f,b],\,[f,c],\,[f,d],\,[f,e]
\end{eqnarray*}

The number of conjugacy classes is $p^{6}+2p^{5}-p^{4}-2p^{3}+p^{2}$ and the
order of the automorphism group is $2(p^{2}-1)^{3}(p^{2}-p)^{3}p^{20}.$

The action of the automorphism group on $G/G^{\prime }$ is given by matrices
in GL$(6,p)$ of the form%
\[
\left( 
\begin{array}{cccccc}
-\alpha & \beta & -\gamma & -\alpha +\delta & \ast & \ast \\ 
\varepsilon & 0 & \zeta & \varepsilon & \ast & \ast \\ 
0 & \eta & 0 & \theta & \ast & \ast \\ 
\alpha & 0 & \gamma & \alpha & \ast & \ast \\ 
0 & 0 & 0 & 0 & \ast & \ast \\ 
0 & 0 & 0 & 0 & \ast & \ast%
\end{array}%
\right) 
\]%
with $(\alpha \zeta -\gamma \varepsilon )(\beta \theta -\delta \eta )\neq 0$
and%
\[
\left( 
\begin{array}{cccccc}
\alpha & -\beta & \gamma & \alpha -\delta & \ast & \ast \\ 
0 & \varepsilon & 0 & \zeta & \ast & \ast \\ 
\eta & 0 & \theta & \eta & \ast & \ast \\ 
0 & \beta & 0 & \delta & \ast & \ast \\ 
0 & 0 & 0 & 0 & \ast & \ast \\ 
0 & 0 & 0 & 0 & \ast & \ast%
\end{array}%
\right) 
\]%
with $(\alpha \theta -\gamma \eta )(\beta \zeta -\delta \varepsilon )\neq 0$.

\bigskip \noindent\textbf{Group 8.6.3}

\begin{eqnarray*}
&&[c,b],\,[d,a],\,[d,b]=[c,a],\,[d,c],\,[e,a],\,[e,b],\,[e,c], \\
&&[e,d],\,[f,a],\,[f,b],\,[f,c],\,[f,d],\,[f,e]
\end{eqnarray*}

The number of conjugacy classes is $p^{6}+p^{5}-p^{3}$ and the order of the
automorphism group is $(p-1)(p^{2}-1)^{2}(p^{2}-p)^{2}p^{24}.$

The action of the automorphism group on $G/G^{\prime }$ is given by matrices
in GL$(6,p)$ of the form%
\[
\left( 
\begin{array}{cccccc}
\alpha & \beta & \ast & \ast & \ast & \ast \\ 
\gamma & \delta & \ast & \ast & \ast & \ast \\ 
0 & 0 & \lambda \delta & -\lambda \gamma & \ast & \ast \\ 
0 & 0 & -\lambda \beta & \lambda \alpha & \ast & \ast \\ 
0 & 0 & 0 & 0 & \ast & \ast \\ 
0 & 0 & 0 & 0 & \ast & \ast%
\end{array}%
\right) 
\]%
with $\lambda (\alpha \delta -\beta \gamma )\neq 0$.

\bigskip \noindent\textbf{Group 8.6.4}

\begin{eqnarray*}
&&[c,b],\,[d,a],\,[d,b]=[c,a],\,[d,c]=[b,a]^{\omega
},\,[e,a],\,[e,b],\,[e,c], \\
&&[e,d],\,[f,a],\,[f,b],\,[f,c],\,[f,d],\,[f,e]
\end{eqnarray*}

The number of conjugacy classes is $p^{6}+p^{4}-p^{2}$ and the order of the
automorphism group is $2(p^{4}-1)(p^{4}-p^{2})(p^{2}-1)(p^{2}-p)p^{20}.$

The action of the automorphism group on $G/G^{\prime }$ is given by matrices
in GL$(6,p)$ of the form%
\[
\left( 
\begin{array}{cccccc}
\alpha & \beta & \gamma & \delta & \ast & \ast \\ 
\varepsilon & \zeta & \eta & \theta & \ast & \ast \\ 
\omega \theta & -\omega \eta & -\zeta & \varepsilon & \ast & \ast \\ 
-\omega \delta & \omega \gamma & \beta & -\alpha & \ast & \ast \\ 
0 & 0 & 0 & 0 & \ast & \ast \\ 
0 & 0 & 0 & 0 & \ast & \ast%
\end{array}%
\right) 
\]%
and 
\[
\left( 
\begin{array}{cccccc}
\alpha & \beta & \gamma & \delta & \ast & \ast \\ 
\varepsilon & \zeta & \eta & \theta & \ast & \ast \\ 
-\omega \theta & \omega \eta & \zeta & -\varepsilon & \ast & \ast \\ 
\omega \delta & -\omega \gamma & -\beta & \alpha & \ast & \ast \\ 
0 & 0 & 0 & 0 & \ast & \ast \\ 
0 & 0 & 0 & 0 & \ast & \ast%
\end{array}%
\right) , 
\]%
where $(\alpha ,\beta ,\gamma ,\delta )$ can be any 4-vector other than zero
($p^{4}-1$ possibilities), and where $(\varepsilon ,\zeta ,\eta ,\theta )$
can be any 4-vector which is \emph{not} in the linear span of $(\alpha
,\beta ,\gamma ,\delta )$ and $(\omega \delta ,-\omega \gamma ,-\beta
,\alpha )$ ($p^{4}-p^{2}$ possibilities).

\bigskip \noindent\textbf{Group 8.6.5}

\begin{eqnarray*}
&&[c,b],\,[d,a],\,[d,b],\,[d,c],\,[e,a],\,[e,b],\,[e,c], \\
&&[e,d]=[b,a],\,[f,a],\,[f,b],\,[f,c],\,[f,d],\,[f,e]
\end{eqnarray*}

The number of conjugacy classes is $p^{6}+p^{5}-p^{3}$ and the order of the
automorphism group is $(p-1)^{3}(p^{2}-1)(p^{2}-p)p^{22}$.

The action of the automorphism group on $G/G^{\prime }$ is given by matrices
in GL$(6,p)$ of the form%
\[
\left( 
\begin{array}{cccccc}
\alpha & \beta & \gamma & \delta & \varepsilon & \ast \\ 
0 & \alpha ^{-1}(\lambda \xi -\mu \nu ) & 0 & 0 & 0 & \ast \\ 
0 & \zeta & \eta & 0 & 0 & \ast \\ 
0 & \alpha ^{-1}(-\delta \mu +\varepsilon \lambda ) & 0 & \lambda & \mu & 
\ast \\ 
0 & \alpha ^{-1}(-\delta \xi +\varepsilon \nu ) & 0 & \nu & \xi & \ast \\ 
0 & 0 & 0 & 0 & 0 & \ast%
\end{array}%
\right) 
\]%
with $\alpha ,\eta ,\lambda \xi -\mu \nu \neq 0$.

\newpage \noindent\textbf{Group 8.6.6}

\begin{eqnarray*}
&&[c,b],\,[d,a],\,[d,b]=[c,a],\,[d,c],\,[e,a],\,[e,b],\,[e,c], \\
&&[e,d]=[b,a],\,[f,a],\,[f,b],\,[f,c],\,[f,d],\,[f,e]
\end{eqnarray*}

The number of conjugacy classes is $p^{6}+p^{4}-p^{2}$ and the order of the
automorphism group is $(p-1)^{2}(p^{2}-1)(p^{2}-p)p^{21}$.

The action of the automorphism group on $G/G^{\prime }$ is given by a
subgroup $H$ of GL$(6,p)$, where the matrices in $H$ have the form%
\[
\left( 
\begin{array}{cccccc}
\alpha & \ast & \ast & \beta & \ast & \ast \\ 
0 & \ast & \ast & 0 & \ast & \ast \\ 
0 & \ast & \ast & 0 & \ast & \ast \\ 
\gamma & \ast & \ast & \delta & \ast & \ast \\ 
0 & \ast & \ast & 0 & \ast & \ast \\ 
0 & 0 & 0 & 0 & 0 & \ast%
\end{array}%
\right) . 
\]%
As we run through the elements of $H$, $\left( 
\begin{array}{cc}
\alpha & \beta \\ 
\gamma & \delta%
\end{array}%
\right) $ takes all values in GL$(2,p)$. If we take generators $\left( 
\begin{array}{cc}
\omega & 0 \\ 
0 & 1%
\end{array}%
\right) $ and $\left( 
\begin{array}{cc}
-1 & 1 \\ 
-1 & 0%
\end{array}%
\right) $ for GL$(2,p)$ then we obtain the following generating matrices for 
$H$:%
\[
\left( 
\begin{array}{cccccc}
\omega & \alpha & \beta & 0 & \gamma & \ast \\ 
0 & \omega \lambda & 0 & 0 & 0 & \ast \\ 
0 & 0 & \lambda & 0 & 0 & \ast \\ 
0 & \omega ^{-1}\gamma & -\omega ^{-1}\alpha & 1 & \delta & \ast \\ 
0 & 0 & 0 & 0 & \omega ^{2}\lambda & \ast \\ 
0 & 0 & 0 & 0 & 0 & \ast%
\end{array}%
\right) \;(\lambda \neq 0) 
\]%
and%
\[
\left( 
\begin{array}{cccccc}
-1 & \alpha & \beta & 1 & \gamma & \ast \\ 
0 & \lambda & 0 & 0 & 2\lambda & \ast \\ 
0 & 0 & 0 & 0 & \lambda & \ast \\ 
-1 & \delta & \beta -\delta & 0 & \delta -\alpha & \ast \\ 
0 & -\lambda & \lambda & 0 & -\lambda & \ast \\ 
0 & 0 & 0 & 0 & 0 & \ast%
\end{array}%
\right) \;(\lambda \neq 0). 
\]

\newpage \noindent\textbf{Group 8.6.7}

\begin{eqnarray*}
&&[c,a],\,[c,b],\,[d,a],\,[d,b],\,[e,a],\,[e,b],\,[e,c], \\
&&[e,d],\,[f,a],\,[f,b],\,[f,c],\,[f,d],\,[f,e]=[b,a]
\end{eqnarray*}

The number of conjugacy classes is $p^{6}+p^{5}-p^{4}+p^{3}-2p+1$ and the
order of the automorphism group is $%
(p^{4}-1)(p^{4}-p^{3})(p^{2}-1)^{2}(p^{2}-p)p^{13}$.

The action of the automorphism group on $G/G^{\prime }$ is given by a
subgroup $H$ of GL$(6,p)$, where the matrices in $H$ have the form%
\[
\left( 
\begin{array}{cccccc}
\ast & \ast & 0 & 0 & \ast & \ast \\ 
\ast & \ast & 0 & 0 & \ast & \ast \\ 
0 & 0 & \ast & \ast & 0 & 0 \\ 
0 & 0 & \ast & \ast & 0 & 0 \\ 
\ast & \ast & 0 & 0 & \ast & \ast \\ 
\ast & \ast & 0 & 0 & \ast & \ast%
\end{array}%
\right) . 
\]%
The first row takes on all possible $p^{4}-1$ non-zero values of the form
shown. Let the first row correspond to an element $g_{1}=a^{\ast }b^{\ast
}e^{\ast }f^{\ast }$. The centralizer of $g_{1}$ has index $p$ in $G$, and
the second row must correspond to an element $g_{2}=a^{\ast }b^{\ast
}e^{\ast }f^{\ast }$ outside the centralizer of $g_{2}$ ($p^{4}-p^{3}$
possibilities). The fifth row must correspond to a non-trivial element $%
g_{5}=a^{\ast }b^{\ast }e^{\ast }f^{\ast }$ which centralizes $g_{1}$ and $%
g_{2}$ (\thinspace $p^{2}-1$ possibilities). The sixth row must correspond
to an element $g_{6}=a^{\ast }b^{\ast }e^{\ast }f^{\ast }$ which centralizes 
$g_{1}$ and $g_{2}$ but does not centralize $g_{5}$ ($p^{2}-p$
possibilities), but we require $[g_{6},g_{5}]=[g_{2},g_{1}]\,$\ and this
reduces the number of choices for $g_{6}$ to $p$. The third and fourth rows
correspond to non-commuting elements of the form $c^{\ast }d^{\ast }$ ($%
(p^{2}-1)(p^{2}-p)$ possibilities).

\bigskip \noindent\textbf{Group 8.6.8}

\begin{eqnarray*}
&&[c,a],\,[c,b],\,[d,a],\,[d,b],\,[e,a],\,[e,b],\,[e,c], \\
&&[e,d],\,[f,a],\,[f,b],\,[f,c],\,[f,d],\,[f,e]=[b,a][d,c]
\end{eqnarray*}

The number of conjugacy classes is $p^{6}+3p^{3}-2p^{2}-3p+2$ and the order
of the automorphism group is $6(p^{2}-1)^{3}(p^{2}-p)p^{14}$.

The action of the automorphism group on $G/G^{\prime }$ is given by a
subgroup $H$ of GL$(6,p)$. The group $H$ has a subgroup $K$ of order $%
(p^{2}-1)^{3}(p^{2}-p)p^{2}$ consisting of matrices of the form%
\[
\left( 
\begin{array}{cccccc}
\ast & \ast & 0 & 0 & 0 & 0 \\ 
\ast & \ast & 0 & 0 & 0 & 0 \\ 
0 & 0 & \ast & \ast & 0 & 0 \\ 
0 & 0 & \ast & \ast & 0 & 0 \\ 
0 & 0 & 0 & 0 & \ast & \ast \\ 
0 & 0 & 0 & 0 & \ast & \ast%
\end{array}%
\right) , 
\]%
where the determinants of the three $2\times 2$ blocks are equal. The three
blocks \textquotedblleft can be permuted around\textquotedblright\ in 6
ways, and the group $H$ is generated by $K$ and the matrices%
\[
\left( 
\begin{array}{cccccc}
0 & 0 & 1 & 0 & 0 & 0 \\ 
0 & 0 & 0 & 1 & 0 & 0 \\ 
1 & 0 & 0 & 0 & 0 & 0 \\ 
0 & 1 & 0 & 0 & 0 & 0 \\ 
0 & 0 & 0 & 0 & 1 & 0 \\ 
0 & 0 & 0 & 0 & 0 & 1%
\end{array}%
\right) \mathrm{ and }\left( 
\begin{array}{cccccc}
0 & 0 & 1 & 0 & 0 & 0 \\ 
0 & 0 & 0 & 1 & 0 & 0 \\ 
0 & 0 & 0 & 0 & -1 & 0 \\ 
0 & 0 & 0 & 0 & 0 & 1 \\ 
-1 & 0 & 0 & 0 & 0 & 0 \\ 
0 & 1 & 0 & 0 & 0 & 0%
\end{array}%
\right) . 
\]

\bigskip \noindent\textbf{Group 8.6.9}

\begin{eqnarray*}
&&[c,b],\,[d,a],\,[d,b]=[c,a],\,[d,c],\,[e,a],\,[e,b],\,[e,c], \\
&&[e,d],\,[f,a],\,[f,b],\,[f,c],\,[f,d],\,[f,e]=[b,a]
\end{eqnarray*}

The number of conjugacy classes is $p^{6}+2p^{3}-p^{2}-2p+1$ and the order
of the automorphism group is $(p^{2}-1)^{2}(p^{2}-p)^{2}p^{15}$.

The action of the automorphism group on $G/G^{\prime }$ is given by the
subgroup of GL$(6,p)$ consisting of matrices%
\[
\left( 
\begin{array}{cccccc}
\alpha & \beta & \varepsilon & \zeta & 0 & 0 \\ 
\gamma & \delta & \eta & \theta & 0 & 0 \\ 
0 & 0 & \lambda \delta & -\lambda \alpha & 0 & 0 \\ 
0 & 0 & -\lambda \beta & \lambda \alpha & 0 & 0 \\ 
0 & 0 & 0 & 0 & \rho & \sigma \\ 
0 & 0 & 0 & 0 & \tau & \varphi%
\end{array}%
\right) 
\]%
with $\alpha \delta -\beta \gamma =\rho \varphi -\sigma \tau \neq 0$, $%
\lambda \neq 0$, $\alpha \eta +\beta \theta =\gamma \varepsilon +\delta
\zeta $.

\newpage \noindent\textbf{Group 8.6.10}

\begin{eqnarray*}
&&[c,b],\,[d,a],\,[d,b]=[c,a],\,[d,c],\,[e,a],\,[e,b],\,[e,c], \\
&&[e,d],\,[f,a],\,[f,b],\,[f,c],\,[f,d],\,[f,e]=[c,a]
\end{eqnarray*}

The number of conjugacy classes is $p^{6}+p^{5}-p^{4}+p^{2}-p$ and the order
of the automorphism group is $(p^{2}-1)^{2}(p^{2}-p)^{2}p^{20}$.

The action of the automorphism group on $G/G^{\prime }$ is given by the
subgroup $H$ of GL$(6,p)$ where the matrices in $H$ have the form%
\[
\left( 
\begin{array}{cccccc}
\ast & \ast & \ast & \ast & \ast & \ast \\ 
\ast & \ast & \ast & \ast & \ast & \ast \\ 
0 & 0 & \ast & \ast & 0 & 0 \\ 
0 & 0 & \ast & \ast & 0 & 0 \\ 
0 & 0 & \ast & \ast & \ast & \ast \\ 
0 & 0 & \ast & \ast & \ast & \ast%
\end{array}%
\right) . 
\]%
There is a subgroup of $H$ of order $p^{8}$ consisting of matrices of the
form%
\[
\left( 
\begin{array}{cccccc}
1 & 0 & \alpha & \beta & \gamma & \delta \\ 
0 & 1 & \varepsilon & \zeta & \eta & \theta \\ 
0 & 0 & 1 & 0 & 0 & 0 \\ 
0 & 0 & 0 & 1 & 0 & 0 \\ 
0 & 0 & \delta & \theta & 1 & 0 \\ 
0 & 0 & -\gamma & -\eta & 0 & 1%
\end{array}%
\right) . 
\]%
The group $H$ is generated by these matrices together with matrices of the
form%
\[
\left( 
\begin{array}{cccccc}
\alpha & \beta & 0 & 0 & 0 & 0 \\ 
\gamma & \delta & 0 & 0 & 0 & 0 \\ 
0 & 0 & \lambda \delta & -\lambda \gamma & 0 & 0 \\ 
0 & 0 & -\lambda \beta & \lambda \alpha & 0 & 0 \\ 
0 & 0 & 0 & 0 & \rho & \sigma \\ 
0 & 0 & 0 & 0 & \tau & \varphi%
\end{array}%
\right) 
\]%
where $\rho \varphi -\sigma \tau =\lambda (\alpha \delta -\beta \gamma )\neq
0$.

\bigskip \noindent\textbf{Group 8.6.11}

\begin{eqnarray*}
&&[c,b],\,[d,a],\,[d,b]=[c,a],\,[d,c]=[b,a]^{\omega
},\,[e,a],\,[e,b],\,[e,c], \\
&&[e,d],\,[f,a],\,[f,b],\,[f,c],\,[f,d],\,[f,e]=[b,a]
\end{eqnarray*}

The number of conjugacy classes is $p^{6}+p^{3}-p$ and the order of the
automorphism group is $2(p^{4}-1)(p^{3}-p^{2})(p^{2}-1)p^{13}$.

The action of the automorphism group on $G/G^{\prime }$ is given by the
matrices in GL$(6,p)$ with the form%
\[
\left( 
\begin{array}{cccccc}
\alpha & \beta & \gamma & \delta & 0 & 0 \\ 
\varepsilon & \zeta & \eta & \theta & 0 & 0 \\ 
\omega \theta & -\omega \eta & -\zeta & \varepsilon & 0 & 0 \\ 
-\omega \delta & \omega \gamma & \beta & -\alpha & 0 & 0 \\ 
0 & 0 & 0 & 0 & \lambda & \mu \\ 
0 & 0 & 0 & 0 & \nu & \xi%
\end{array}%
\right) \mathrm{ and }\left( 
\begin{array}{cccccc}
\alpha & \beta & \gamma & \delta & 0 & 0 \\ 
\varepsilon & \zeta & \eta & \theta & 0 & 0 \\ 
-\omega \theta & \omega \eta & \zeta & -\varepsilon & 0 & 0 \\ 
\omega \delta & -\omega \gamma & -\beta & \alpha & 0 & 0 \\ 
0 & 0 & 0 & 0 & \lambda & \mu \\ 
0 & 0 & 0 & 0 & \nu & \xi%
\end{array}%
\right) 
\]%
where $\alpha \eta +\beta \theta =\gamma \varepsilon +\delta \zeta $ and $%
\lambda \xi -\mu \nu =\alpha \zeta -\beta \varepsilon +\omega (\gamma \theta
-\delta \eta )\neq 0$. Note that there are $p^{4}-1$ choices for row 1, and
that once row 1 has been fixed there are $p^{3}-p^{2}$ choices for row 2.

\bigskip \noindent\textbf{Group 8.6.12}

\begin{eqnarray*}
&&[c,b],\,[d,a],\,[d,b],\,[d,c],\,[e,a],\,[e,b],\,[e,c], \\
&&[e,d]=[b,a],\,[f,a],\,[f,b],\,[f,c],\,[f,d]=[c,a],\,[f,e]
\end{eqnarray*}

The number of conjugacy classes is $p^{6}+p^{4}-p^{2}$ and the order of the
automorphism group is $(p^{2}-1)^{2}(p^{2}-p)^{2}p^{18}$.

The action of the automorphism group on $G/G^{\prime }$ is given by a
subgroup $H$ of GL$(6,p)$, where the elements of have the form%
\[
\left( 
\begin{array}{cccccc}
\alpha  & \ast  & \ast  & \beta  & \ast  & \ast  \\ 
0 & \ast  & \ast  & 0 & \ast  & \ast  \\ 
0 & \ast  & \ast  & 0 & \ast  & \ast  \\ 
\gamma  & \ast  & \ast  & \delta  & \ast  & \ast  \\ 
0 & \ast  & \ast  & 0 & \ast  & \ast  \\ 
0 & \ast  & \ast  & 0 & \ast  & \ast 
\end{array}%
\right) 
\]%
where $\left( 
\begin{array}{cc}
\alpha  & \beta  \\ 
\gamma  & \delta 
\end{array}%
\right) $ takes all values in GL$(2,p)$. There is a subgroup of $H$ of order 
$p^{6}$ consisting of all elements of the form%
\[
\left( 
\begin{array}{cccccc}
1 & \varepsilon  & \zeta  & 0 & \eta  & \theta  \\ 
0 & 1 & 0 & 0 & 0 & 0 \\ 
0 & 0 & 1 & 0 & 0 & 0 \\ 
0 & \eta  & \theta  & 1 & \lambda  & \mu  \\ 
0 & 0 & 0 & 0 & 1 & 0 \\ 
0 & 0 & 0 & 0 & 0 & 1%
\end{array}%
\right) ,
\]%
and $H$ is generated by this subgroup together with elements%
\[
\left( 
\begin{array}{cccccc}
\alpha  & 0 & 0 & \beta  & 0 & 0 \\ 
0 & \delta \lambda  & \delta \mu  & 0 & -\gamma \lambda  & -\gamma \mu  \\ 
0 & \delta \nu  & \delta \xi  & 0 & -\gamma \nu  & -\gamma \xi  \\ 
\gamma  & 0 & 0 & \delta  & 0 & 0 \\ 
0 & -\beta \lambda  & -\beta \mu  & 0 & \alpha \lambda  & \alpha \mu  \\ 
0 & -\beta \nu  & -\beta \xi  & 0 & \alpha \nu  & \alpha \xi 
\end{array}%
\right) 
\]%
with $(\alpha \delta -\beta \gamma )(\lambda \xi -\mu \nu )\neq 0$.

\bigskip \noindent\textbf{Group 8.6.13}

\begin{eqnarray*}
&&[b,a],\,[d,a],\,[e,a][c,a],\,[f,a],\,[c,b],\,[d,b]=[c,a],\,[e,b], \\
&&[f,b][c,a]^{2},\,[d,c],\,[e,c],\,[e,d]=[f,c],\,[f,d],\,[f,e]=[c,a][f,c]
\end{eqnarray*}

The number of conjugacy classes is $p^{6}+p^{3}-p$ and the order of the
automorphism group is $(p-1)(p^{2}-1)(p^{2}-p)p^{19}$.

The action of the automorphism group on $G/G^{\prime }$ is given by a
subgroup $H=KL$ of GL$(6,p)$, where $K$ is a subgroup of $H$ and $L$ is a
normal subgroup of $H$. The subgroup $K$ is the set of all matrices of the
form%
\[
\left( 
\begin{array}{cccccc}
\alpha & \beta & 0 & 0 & 0 & 0 \\ 
\gamma & \delta & 0 & 0 & 0 & 0 \\ 
0 & 0 & \lambda ^{2}\delta & -\lambda ^{2}\gamma & 0 & 0 \\ 
0 & 0 & -\lambda ^{2}\beta & \lambda ^{2}\alpha & 0 & 0 \\ 
0 & 0 & (\lambda -\lambda ^{2})\delta & (2\lambda +\lambda ^{2})\gamma & 
\lambda \delta & \lambda \gamma \\ 
0 & 0 & (\lambda +2\lambda ^{2})\beta & 2(\lambda -\lambda ^{2})\alpha & 
\lambda \beta & \lambda \alpha%
\end{array}%
\right) 
\]%
with $\lambda (\alpha \delta -\beta \gamma )\neq 0$, and $L$ is the set of
all matrices of the form$\allowbreak $

\[
\left( 
\begin{array}{cccccc}
1 & 0 & 0 & 0 & 0 & 0 \\ 
0 & 1 & 0 & 0 & 0 & 0 \\ 
\alpha & -2\beta ^{2}+2\beta +\varepsilon -\zeta +\gamma \delta & \beta +1 & 
\gamma & \beta & \frac{1}{2}\gamma \\ 
-\beta ^{2}+2\beta +\varepsilon -\zeta +\frac{1}{2}\gamma \delta & \frac{1}{2%
}\delta -\frac{1}{2}\eta & \delta & 2\beta +1 & \delta & \beta \\ 
-\alpha -\frac{1}{2}\gamma & \varepsilon & -\beta & -\gamma & 1-\beta & -%
\frac{1}{2}\gamma \\ 
\zeta & \eta & -2\delta & -4\beta & -2\delta & 1-2\beta%
\end{array}%
\right) . 
\]

\newpage \noindent\textbf{Group 8.6.14}

\begin{eqnarray*}
&&[c,b],\,[d,a],\,[d,b]=[c,a],\,[d,c],\,[e,a],\,[e,b],\,[e,c], \\
&&[e,d]=[b,a],\,[f,a],\,[f,b]=[c,a]^{m},\,[f,c]=[b,a],\,[f,d],\,[f,e]=[c,a],
\end{eqnarray*}%
where $x^{3}-mx+1$ is irreducible over GF$(p)$. (Different choices of $m$
give isomorphic groups.)

The number of conjugacy classes is $p^{6}+p^{2}-1$ and the order of the
automorphism group is $3(p^{6}-1)(p-1)p^{15}$.

Note that since $x^{3}-mx+1$ is irreducible over GF$(p)$ its discriminant $%
4m^{3}-27$ must be a square, and we let $u^{2}=4m^{3}-27$.

The action of the automorphism group on $G/G^{\prime }$ is given by a
subgroup $H$ GL$(6,p)$. The first rows of the matrices in $H$ are completely
arbitrary, except that they must be non-zero. The subgroup $K$ of $H$
consisting of those matrices in $H$ with first row $(1,0,0,0,0,0)$ has order 
$3(p-1)p^{3}$ and is generated by the following matrices%
\[
\left( 
\begin{array}{cccccc}
1 & 0 & 0 & 0 & 0 & 0 \\ 
0 & \lambda & 0 & 0 & 0 & 0 \\ 
0 & 0 & \lambda & 0 & 0 & 0 \\ 
0 & 0 & 0 & 1 & 0 & 0 \\ 
0 & 0 & 0 & 0 & \lambda & 0 \\ 
0 & 0 & 0 & 0 & 0 & 1%
\end{array}%
\right) ,\;\left( 
\begin{array}{cccccc}
1 & 0 & 0 & 0 & 0 & 0 \\ 
-m\gamma -\beta & 1 & 0 & -\gamma & 0 & -\alpha \\ 
\alpha & 0 & 1 & \beta & 0 & \gamma \\ 
0 & 0 & 0 & 1 & 0 & 0 \\ 
-\gamma & 0 & 0 & m\beta -\alpha & 1 & -\beta \\ 
0 & 0 & 0 & 0 & 0 & 1%
\end{array}%
\right) 
\]

and%
\[
A=\left( 
\begin{array}{cccccc}
1 & 0 & 0 & 0 & 0 & 0 \\ 
0 & \frac{1}{2u}(9-u) & \frac{m^{2}}{u} & 0 & \frac{m}{u}(u-3) & 0 \\ 
0 & -\frac{3m}{u} & -\frac{1}{2u}(9+u) & 0 & \frac{2m^{2}}{u} & 0 \\ 
\frac{m}{2u}(3-u) & 0 & 0 & \frac{1}{2u}(9-u) & 0 & -\frac{m^{2}}{u} \\ 
0 & 0 & 0 & 0 & 1 & 0 \\ 
\frac{m^{2}}{u} & 0 & 0 & \frac{3m}{u} & 0 & -\frac{1}{2u}(9+u)%
\end{array}%
\right) . 
\]%
(Here $A^{3}=I$.)

The group $H$ is generated by the matrices above, together with $%
(p^{6}-1)(p^{4}-p^{3})$ matrices of the form%
\[
B=\left( 
\begin{array}{cccccc}
\alpha & \beta & \gamma & \delta & \varepsilon & \zeta \\ 
\mu +m\rho & \lambda & m\xi +\sigma & \rho & -\xi & \nu \\ 
-\nu & \xi & \lambda & -\mu & \sigma & -\rho \\ 
\zeta & \varepsilon & -\beta & \alpha -m\delta & -\gamma -m\varepsilon & 
\delta \\ 
\rho & -\sigma & \xi & \nu -m\mu & \lambda +m\sigma & \mu \\ 
\delta +m\zeta & -\gamma & -m\beta -\varepsilon & \zeta & \beta & \alpha%
\end{array}%
\right) . 
\]

Note that the first row of $B$ determines rows 4 and 6, and that the third
of $B$ determines rows 2 and 5. The first row can take any non-zero value ($%
p^{6}-1$ possibilities). Once the first row of $B$ is fixed the entries $%
\lambda ,\mu ,\nu ,\xi \rho ,\sigma $ from row 3 are required to satisfy two
homogeneous linear equations with coefficients determined by the first row.
These two linear equations correspond to the requirement that the group
elements corresponding to rows 3 and 4 commute. There are $p^{4}$ solutions
to these two equations, but we also require that the third row lie outside
the subspace spanned by rows 1, 4 and 6. Any row 3 lying in this subspace
automatically satisfies these two linear equations, since the group elements
corresponding to rows 1, 4 and 6 are guaranteed to commute by the choice of
rows 4 and 6. So, once row 1 of $B$ is fixed, there are $p^{4}-p^{3}$
choices for row 3 of $B$, and then all the rows of $B$ are determined.

It follows from all this that $H$ has order $3(p^{6}-1)(p-1)p^{3}$.

If we want to obtain explicit generators for $H$ then we only need to find
one possible third row of $B$ for any given first row. Then these particular
choices of $B$, together with the generators of $K$ will generate $H$. If
one of $\alpha ,\delta ,\zeta $ is non-zero, then we can take the third row
of $B$ to be%
\[
(0,\delta ^{2}+m\zeta \delta -\alpha \zeta ,-\alpha ^{2}+m\delta \alpha
+\delta \zeta ,0,\zeta ^{2}-\alpha \delta ,0), 
\]%
and if one of $\beta ,\gamma ,\varepsilon $ is non-zero, then we can take
the third row of $B$ to be%
\[
(\gamma ^{2}+m\varepsilon \gamma -\beta \varepsilon ,0,0,\gamma \varepsilon
+\beta ^{2},0,\varepsilon ^{2}+m\beta \varepsilon +\beta \gamma ). 
\]%
Experimentally, it seems that that $H$ is generated by the generators of $K$
together with the single matrix%
\[
\left( 
\begin{array}{cccccc}
0 & 1 & 0 & 0 & 0 & 0 \\ 
-1 & 0 & 0 & 0 & 0 & 0 \\ 
0 & 0 & 0 & 1 & 0 & 0 \\ 
0 & 0 & -1 & 0 & 0 & 0 \\ 
0 & 0 & 0 & m & 0 & -1 \\ 
0 & 0 & -m & 0 & 1 & 0%
\end{array}%
\right) . 
\]

\subsection{Seven generator groups}

\noindent\textbf{Group 8.7.1}

\[
\langle a,b\rangle \times \langle c\rangle \times \langle d\rangle \times
\langle e\rangle \times \langle f\rangle \times \langle g\rangle . 
\]

The number of conjugacy classes is $p^{7}+p^{6}-p^{5}$, and the automorphism
group has order $%
(p^{8}-p^{6})(p^{8}-p^{7})(p^{6}-p)(p^{6}-p^{2})(p^{6}-p^{3})(p^{6}-p^{4})(p^{6}-p^{5}) 
$.

\bigskip \noindent\textbf{Group 8.7.2}

\[
\langle a,b\rangle \times _{\lbrack b,a]=[d,c]}\langle c,d\rangle \times
\langle e\rangle \times \langle f\rangle \times \langle g\rangle . 
\]

The number of conjugacy classes is $p^{7}+p^{4}-p^{3}$, and the automorphism
group has order $%
(p^{8}-p^{4})(p^{8}-p^{7})(p^{6}-p^{4})p^{5}(p^{4}-p)(p^{4}-p^{2})(p^{4}-p^{3}) 
$.

\bigskip \noindent\textbf{Group 8.7.3}

\[
\langle a,b\rangle \times _{\lbrack b,a]=[d,c]=[f,e]}\langle c,d\rangle
\times _{\lbrack b,a]=[d,c]=[f,e]}\langle e,f\rangle \times \langle g\rangle
. 
\]

The number of conjugacy classes is $p^{7}+p^{2}-p$, and the automorphism
group has order $%
(p^{8}-p^{2})(p^{8}-p^{7})(p^{6}-p^{2})p^{5}(p^{4}-p^{2})p^{3}(p^{2}-p)$.%

\end{document}